\documentclass[11pt,reqno]{amsart}
\usepackage{amsmath,amsthm,amssymb,mathrsfs,stmaryrd,color}
\usepackage[all]{xy}
\usepackage{url}
\usepackage{geometry}
\usepackage[utf8]{inputenc}
\usepackage[T1]{fontenc}

\usepackage{graphicx}

\usepackage{enumitem}
\setlist[enumerate]{itemsep=2pt,parsep=2pt,before={\parskip=2pt}}

\usepackage[colorlinks=true,hyperindex, linkcolor=magenta, pagebackref=false, citecolor=cyan]{hyperref}

\usepackage{relsize}
\usepackage[bbgreekl]{mathbbol}
\usepackage{amsfonts}
\DeclareSymbolFontAlphabet{\mathbb}{AMSb}
\DeclareSymbolFontAlphabet{\mathbbl}{bbold}
\newcommand{\Prism}{{\mathlarger{\mathbbl{\Delta}}}}

\usepackage{slashed}
\DeclareMathOperator{\DHod}{\slashed{D}}

\begin{document}

\newtheorem{theorem}{Theorem}[section]
\newtheorem*{theorem*}{Theorem}
\newtheorem*{condition*}{Condition}
\newtheorem*{definition*}{Definition}
\newtheorem{proposition}[theorem]{Proposition}
\newtheorem{lemma}[theorem]{Lemma}
\newtheorem{corollary}[theorem]{Corollary}
\newtheorem{claim}[theorem]{Claim}
\newtheorem{claimex}{Claim}[theorem]
\newtheorem{conjecture}[theorem]{Conjecture}
\newtheorem{etheorem}[theorem]{Expected Theorem}

\theoremstyle{definition}
\newtheorem{definition}[theorem]{Definition}
\newtheorem{question}[theorem]{Question}
\newtheorem{observation}[theorem]{Observation}
\newtheorem{remark}[theorem]{Remark}
\newtheorem{construction}[theorem]{Construction}
\newtheorem{example}[theorem]{Example}
\newtheorem{condition}[theorem]{Condition}
\newtheorem{warning}[theorem]{Warning}
\newtheorem{notation}[theorem]{Notation}
\newtheorem*{notation*}{Notation}
\newtheorem{goal}[theorem]{Goal}
\newtheorem{problem}[theorem]{Problem}
\newtheorem{guess}[theorem]{Guess}

\title{Algebraic geometry in mixed characteristic}
\author{Bhargav Bhatt}

\begin{abstract}
Fix a prime number $p$. We report on some recent developments in algebraic geometry (broadly construed) over $p$-adically complete commutative rings. These developments include foundational advances within the subject as well as external applications. 
\end{abstract}

\maketitle
\setcounter{tocdepth}{1}
\tableofcontents


\section*{Introduction}

Fix a prime number $p$. The last decade has witnessed multiple conceptual advances in algebraic geometry over mixed characteristic rings (which, in this article, we take to mean commutative rings that are $p$-adically complete). These advances have led to the resolution of longstanding questions in different areas of mathematics where $p$-adic completions appear. Moreover, entirely new and fascinating domains of inquiry have been uncovered.  The goal of this survey is to discuss some of these developments, especially in topics close to the author's expertise.

A highlight of the last decade of activity in the area has been its seat as an exchange of ideas across different fields of mathematics.  For instance, a central topic of this survey is prismatic cohomology, which is a new cohomology theory for mixed characteristic rings (\S \ref{sec:RelativePrismatic} and \S \ref{sec:APC}); its discovery was inspired in part by calculations in homotopy theory and in part by developments in Galois representation theory (Remark~\ref{rmk:OriginKMot}). Prismatic cohomology in turn played a prominent role in the proof of a mixed characteristic analog of the Kodaira vanishing theorem (Theorem~\ref{thm:KV}), which then helped develop the minimal model program in the birational geometry of arithmetic threefolds (Theorem~\ref{MMP}). In the reverse direction, an important flatness lemma discovered in the solution \cite{AndreDSC} of a longstanding question in commutative algebra facilitated, via prismatic cohomology again, the proof of an analog of Bott's vanishing theorem for algebraic $K$-theory (Theorem~\ref{OddVan}). The author hopes this survey can convey some of the excitement surrounding this interplay of ideas.

We emphasize right away that the topics covered are chosen somewhat idiosyncratically, and we have not attempted to be comprehensive even in the topics we do cover; to make partial amends, a number of references have been included to help the reader navigate the subject. Moreover, the level of the exposition is uneven across sections; for instance, we have taken a macroscopic view of topics that are reasonably well documented elsewhere, but have gone into more detail and depth while covering very recent ideas that seem promising. 

This survey is organized as follows. In \S \ref{sec:RelativePrismatic}, we discuss relative prismatic cohomology and related developments. The absolute version of this story, which is comparatively new, is the subject of \S \ref{sec:APC}. We then present applications, covering algebraic $K$-theory in \S \ref{sec:K} and commutative algebra and birational geometry in \S \ref{sec:CABG}.  We end in \S \ref{sec:RH} with some relatively recent work on $p$-adic Riemann--Hilbert problems and their algebro-geometric implications.

All rings that appear are assumed commutative unless otherwise specified.

\section{Prisms and relative prismatic cohomology}
\label{sec:RelativePrismatic}
In the last century, especially following the work of Grothendieck, cohomology theories have emerged as extremely important tools in algebraic geometry and number theory: they lie at the heart of some of the deepest theorems and conjectures in both subjects. For example, classical Hodge theory, which studies the singular cohomology with real/complex coefficients for complex varieties, is a central topic in modern algebraic geometry, with applications throughout the subject and beyond. Likewise, $p$-adic Hodge theory, which studies the $p$-adic cohomology of $p$-adic varieties, is an equally fundamental notion in arithmetic geometry:  it provides one of the best known tools for understanding Galois representations of the absolute Galois group of $\mathbf{Q}$. Moreover, unlike in the complex setting, there is a large number of cohomology theories in the $p$-adic world:  \'etale, de Rham, Hodge, crystalline, de Rham--Witt, etc. In this section, we will report on work from the last few years dedicated to finding 
an organizational framework to better understand $p$-adic cohomology theories in $p$-adic arithmetic geometry, especially their relationships with each other.

\begin{remark}[Why do derived objects appear repeatedly?]
\label{rmk:derived}
Before embarking on our journey, let us explain one reason derived notions (i.e., those with a homological/homotopical flavor) often appear in recent work in the area and consequently also in our exposition.

In classical algebraic geometry, the fundamental objects are smooth algebraic varieties over an algebraically closed field. Similarly, in mixed characteristic algebraic geometry, the basic geometric objects are ($p$-adic formal\footnote{A $p$-adic formal scheme is a formal scheme whose affine opens are given by formal spectra of $p$-adically complete rings equipped with the $p$-adic topology}) smooth schemes over the ring of integers $\mathcal{O}_C$ of a complete algebraically closed nonarchimedean field $C/\mathbf{Q}_p$. In particular, unlike the classical setting, the rings of functions that appear in mixed characteristic algebraic geometry are often not noetherian: indeed, the ring $\mathcal{O}_C$ is a non-noetherian valuation ring as its value group is divisible. Replacing $\mathcal{O}_C$ with a discrete valuation ring like $\mathbf{Z}_p$, while quite tempting and important for applications, leads to arithmetic subtleties that one would like to avoid, at least at first pass, in a purely geometric study. Even more exotic non-noetherian rings are critical to several recent innovations in the area, such as perfectoid geometry \cite{ScholzePerfectoid,KedlayaLiu1}, descent techniques for extremely fine Grothendieck topologies such as the pro-\'etale, quasi-syntomic, $v$ and arc topologies \cite{BhattScholzeProetale,ScholzeDiamonds,BhattScholzeWitt,BhattMathewArc,BMS2}, the theory of $\delta$-rings \cite{Joyal}, etc.

In the non-classical situations described above, derived notions often have better stability properties than their classical counterparts. For instance, given a commutative ring $R$ with a finitely generated ideal $I$, the category of derived $I$-complete $R$-modules forms an abelian subcategory (e.g., it is closed under kernels and cokernels) of the category of all $R$-modules, unlike the subcategory of classically $I$-adically complete $R$-modules; moreover, the assignment carrying $R$ to the $\infty$-category $\mathcal{D}_{I-\text{comp}}(R)$ of derived $I$-complete $R$-complexes forms a stack for the flat topology (or even a suitably defined $I$-completely flat topology), unlike the corresponding assignment at the triangulated category level. For such reasons, the language of higher category theory and derived algebraic geometry \cite{LurieHTT,LurieHA,LurieSAG,ToenVezz1,ToenVezz2} has played an important role in the developments discussed in this paper.
\end{remark}

The work described in this section began with the goal to enhance Fontaine's perspective \cite{FontaineAnnals} on $p$-adic Hodge theory to work well with integral coefficients. A concrete goal was to understand how the torsion in the $\mathbf{Z}_p$-\'etale cohomology of the geometric generic fibre of a smooth projective scheme over a mixed characteristic discrete valuation ring relates to the torsion in the crystalline cohomology of its special fibre; this question was already stressed by Grothendieck in his Algerian letter to Deligne in $1965$ (see \cite{BreuilMessingTorsion} for a survey on the status of this question $20$ years ago, and \S 2 in {\em loc.\ cit.} for history).  After initial attempts \cite{BMS1,BMS2} that worked well in important examples, a satisfactory theory was found via the notion of a prism, recalled next. The definition relies on the notion of a $\delta$-ring, which is roughly a commutative ring $A$ equipped with a map $\varphi:A \to A$ lifting the Frobenius endomorphism $f \mapsto f^p$ of $A/p$ (interpreted in a derived sense when $A$ has $p$-torsion); see \cite{Joyal,BuiumDerivation}. The importance of this notion in arithmetic geometry has long been stressed by Borger, see \cite{BorgerWitt}.

\begin{definition}[Prisms, \cite{BhattScholzePrisms}]
A {\em prism} is a pair $(A,I)$, where $A$ is a $\delta$-ring and $I \subset A$ is an invertible ideal such that $A$ is derived $(p,I)$-complete and $p \in (I,\varphi(I))$. Write $\overline{A} := A/I$.
\end{definition}

In practice, we restrict to {\em bounded} prisms, i.e., those prisms $(A,I)$ where the $p$-power torsion in $\overline{A}$ is annihilated by $p^n$ for some $n \geq 0$; this restriction allows us to avoid certain derived technicalities without sacrificing the key examples. Two important examples are discussed next; see Remark~\ref{rmk:qdR} for another key example.

\begin{example}[Crystalline prisms] 
\label{CrysPrism}
If $A$ is any $p$-complete $p$-torsionfree $\delta$-ring, then $(A,(p))$ is a bounded prism. For instance, given a reduced $\mathbf{F}_p$-algebra $R$, we could take $A=W(R)$ to be the ring of $p$-typical Witt vectors of ring with its natural Frobenius lift.
\end{example}

\begin{example}[Perfect prisms]
\label{PerfPrism}
 A prism $(A,I)$ is called {\em perfect} if the Frobenius $\varphi:A \to A$ is an isomorphism; any such prism is bounded and the ring $\overline{A}$ is perfectoid as in \cite{GabberRameroFART,BMS1}. In fact, the construction $(A,I) \mapsto \overline{A}$ yields an equivalence of categories between perfect prisms and perfectoid rings; thus, the notion of a prism may be viewed as a ``deperfection'' of the notion of a perfectoid ring.  An important example is the perfect prism $(A,I)$ corresponding to the perfectoid ring $\overline{A} = \mathcal{O}_C$, where $C/\mathbf{Q}_p$ is a complete and algebraically closed extension; we call this {\em a Fontaine prism}, in homage to its discovery \cite[\S 5]{FontainepdivBook}.
\end{example}

Given a bounded prism $(A,I)$ as well as an $\overline{A}$-scheme $X$, the following key definition allows us to extract an $A$-linear cohomology theory for $X$.

\begin{definition}[The relative prismatic site]
\label{def:RelPrismSite}
Fix a bounded prism $(A,I)$ and a $p$-adic formal $\overline{A}$-scheme $X$. The relative prismatic site $(X/A)_\Prism$ is the category of all bounded prisms $(B,J)$ over $(A,I)$ equipped with an $\overline{A}$-map $\mathrm{Spf}(B/J) \to X$, topologized via the flat topology; write $\mathcal{O}_\Prism$, $I_\Prism$ and $\overline{\mathcal{O}}_\Prism$ for the sheaves obtained by remembering $B$, $J$ or $B/J$ respectively, so there is a natural $\varphi$-action on $\mathcal{O}_\Prism$ and $\overline{\mathcal{O}}_\Prism = \mathcal{O}_\Prism/I_\Prism$. Write $R\Gamma_\Prism(X/A) := R\Gamma( (X/A)_\Prism, \mathcal{O}_\Prism) \ \rotatebox[origin=c]{-270}{$\circlearrowright$} \ \varphi$ for the resulting cohomology theory.
\end{definition}

The main comparison theorems for $R\Gamma_\Prism(X/A)$ are informally summarized next:

\begin{theorem}[Relative prismatic cohomology, \cite{BMS1,BhattScholzePrisms}]
\label{PrismaticCoh}
Fix a bounded prism $(A,I)$ and let $X$ be a smooth $p$-adic formal $\overline{A}$-scheme.  The relative prismatic cohomology theory $R\Gamma_\Prism(X/A)  \ \rotatebox[origin=c]{-270}{$\circlearrowright$} \ \varphi$ recovers the standard integral $p$-adic cohomology theories for $X/A$ with their extra structures (e.g., \'etale, de Rham, Hodge, crystalline, de Rham--Witt) via a specialization procedure, thereby giving new relationships between them. 
\end{theorem}

For instance, if $(A,I)$ is crystalline, Theorem~\ref{PrismaticCoh} leads to a canonical Frobenius descent of crystalline cohomology \cite{BerthelotCrys}; this descent was previously observed on cohomology groups in \cite{OgusVologodsky,ShihoOV}. On the other hand, the $A_{\inf}$-cohomology of \cite{BMS1} is recovered by specializing to a Fontaine prism; we refer to the surveys  \cite{MorrowNotesAinf,BhattSpecializing} for more precise assertions (with pictures!) about the comparisons in this case. An early concrete application of the latter was the following result relating \'etale and de Rham cohomology integrally; via classical comparisons, this gives a new technique to bound the $p$-torsion in singular cohomology of complex algebraic varieties via the geometry of their mod $p$ reductions.

\begin{corollary}[Torsion inequality, \cite{BMS1}]
Let $C/\mathbf{Q}_p$ be a complete and algebraically closed field with residue field $k$ (e.g., we may take $C=\mathbf{C}_p=\widehat{\overline{\mathbf{Q}_p}}$, so $k=\overline{\mathbf{F}_p}$). Let $X/\mathcal{O}_C$ be a proper smooth $p$-adic formal scheme. Then 
\[ \dim_{\mathbf{F}_p} H^i_{et}(X_C,\mathbf{F}_p) \leq \dim_k H^i_{\mathrm{dR}}(X_k) \quad \text{for all } i \geq 0.\]
More generally, a similar inequality bounds the length of the torsion subgroup of $H^i_{et}(X_C,\mathbf{Z}_p)$ in terms of that of $H^i_{\text{crys}}(X_k)$. In particular, if the latter is torsionfree, so is the former.
\end{corollary}

Since its discovery, the prismatic theory in \cite{BMS1,BhattScholzePrisms} has found several applications, some of which are discussed below and elsewhere in this paper. Other results featuring this theory include: Hodge theory of classifying spaces of reductive groups \cite{KPTotaro, BhattLiTotaro}, vanishing theorems for the cohomology of the moduli space of curves with level structures \cite{ReineckeMg}, essential dimension calculations \cite{FKWEssentialDimAbVar}, Poincar\'e duality for $\mathbf{Z}/p^n$-coefficients in rigid geometry \cite{ZavyalovPD}, calculation of the $\mathbf{Z}_p$-cohomology of Drinfeld's $p$-adic symmetric spaces \cite{CDNDrinfeld}, a fairly optimal form of Dieudonn\'e theory in mixed characteristic \cite{ALBPrismatic}, better understanding of the moduli stacks of Breuil-Kisin-Fargues modules \cite{EmertonGeeStack}, and several improvements to known comparisons in integral $p$-adic Hodge theory \cite{MinpHT, LipHT, LiLiuCompare}.

\begin{remark}[Rational comparison theorems]
Specializing part of Theorem~\ref{PrismaticCoh} to a Fontaine prism $(A,I)$ gives a generalization of Fontaine's crystalline comparison conjecture $C_{\mathrm{crys}}$ to proper smooth formal schemes $X/\overline{A}$; variants of both this result and its proof have a long history in $p$-adic Hodge theory, including \cite{FontaineMessing,BlochKatoComp,FaltingspHT,FaltingsVeryRamified,TsujiComp,NiziolComp,BeilinsonComp,ColmezNiziolComp}.
\end{remark}

\begin{remark}[$q$-de Rham cohomology, \cite{BhattScholzePrisms}]
\label{rmk:qdR}
Given a smooth $\mathbf{Z}$-algebra $R$ equipped with a choice of \'etale co-ordinates (which we call a {\em framing} and indicate by $\square$), one can define a complex $q\Omega^\bullet_{(R,\square)}$ of $\mathbf{Z}\llbracket q-1 \rrbracket$-modules by $q$-deforming the differential of the de Rham complex $\Omega^*_{R/\mathbf{Z}}$ (see \cite{Aomotoq,ScholzeqdR}); this complex strongly depends on the framing $\square$. Nevertheless, motivated by some local calculations from \cite{BMS1}, Scholze had conjectured in \cite{ScholzeqdR} that $q\Omega^\bullet_{(R,\square)}$ is independent of the framing $\square$ up to canonical quasi-isomorphism. This conjecture was deduced from the existence of prismatic cohomology in \cite{BhattScholzePrisms}, as explained next; prior partial progress was made by Pridham \cite{Pridhamq}, also using $\delta$-rings.

By a patching procedure, Scholze's co-ordinate independence conjecture reduces to its analog when all objects are $p$-completed. The latter follows from the existence of prismatic cohomology relative to the $q$-de Rham prism $(A,I) := (\mathbf{Z}_p\llbracket q-1\rrbracket, (\frac{q^p-1}{q-1}))$ where  $\varphi(q) = q^p$: given a formally smooth $p$-complete $\mathbf{Z}_p$-algebra $R$ equipped with a framing $\square$ as before, the relative prismatic complex $R\Gamma_{\Prism}(\mathrm{Spf}(R \otimes_{\mathbf{Z}_p} \overline{A})/A)$ (which is visibly independent of the framing $\square$) is naturally quasi-isomorphic to the $q$-dR complex $q\Omega^\bullet_{(R,\square)}$. 

The preceding perspective on $q$-de Rham cohomology also yields a formalism for more systematically discussing related notions such as, e.g., Gauss-Manin $q$-connections; we refer to \cite{Chatzistamatiouq,GLSQq1,GLSQq2,MorrowTsuji} for more on these and related developments.
\end{remark}

\begin{remark}[Perfections in mixed characteristic, \cite{BhattScholzePrisms}]
\label{rmk:Perfections}
The theory of perfectoid rings can be regarded as a mixed characteristic analog of the theory of perfect $\mathbf{F}_p$-algebras i.e., $\mathbf{F}_p$-algebras where the Frobenius is bijective. The utility of this analogy is enhanced by Theorem~\ref{PrismaticCoh}: by attaching objects with Frobenius actions to rings in mixed characteristic, this result yields a notion of ``perfectoidization'' for a large class of mixed characteristic rings. Indeed, given any perfect prism $(A,I)$ and an $\overline{A}$-algebra $R$, one can naturally construct a ``(derived) perfectoidization'' $R \to R_{\mathrm{perfd}}$ with excellent formal properties. For instance, if $R$ is integral over $\overline{A}$, then $R \to R_{\mathrm{perfd}}$ is in fact the universal map to a perfectoid ring.   This construction has several applications. For instance, \cite{BhattScholzePrisms} uses these to prove an optimal generalization of the Faltings' almost purity theorem (extending versions from \cite{FaltingspHT,FaltingsAEE,ScholzePerfectoid,KedlayaLiu1,AndreAbhyankar}) as well as the result that ``Zariski closed = strongly Zariski closed'' for affinoid perfectoid spaces; the latter plays an important role in aspects of \cite{FarguesScholze}. The perfectoidization functor also powers the construction of the $p$-adic Riemann--Hilbert functor in \S \ref{sec:RH}.
\end{remark}

\begin{remark}[A new perspective on de Rham--Witt complexes, \cite{BLMDRW}]
\label{rmk:DRW}
The de Rham-Witt complex of Bloch--Deligne--Illusie \cite{BlochKCrys,IllusiedRW} is a fundamental object in characteristic $p$ algebraic geometry with applications transcending  algebra (e.g., \cite{HesselholtMadsen}). Its construction traditionally relied on somewhat laborious calculations. The paper \cite{BLMDRW}, inspired by structures on relative prismatic cohomology, offered a new homological perspective on this object. 

To explain this, we first recall the isogeny theorem for prismatic cohomology. In the setup of Theorem~\ref{PrismaticCoh}, when $X=\mathrm{Spf}(R)$ is affine, one often writes $\Prism_{R/A} = R\Gamma_\Prism(X/A)$, regarded as an object of the derived category of $A$. There is then a natural quasi-isomorphism 
\begin{equation}
\tag{Isog}
\label{eq:Leta}
 \widetilde{\varphi}_{R/A}:\varphi^* \Prism_{R/A} \simeq L\eta_I \Prism_{R/A}
 \end{equation}
induced by the relative Frobenius, where $L\eta_I$ is a variant of the Berthelot-Ogus-Deligne d\'ecalage functor (see \cite[\S 6]{BMS1}); the isomorphism $\widetilde{\varphi}_{R/A}$, which is a prismatic avatar of the Berthelot-Ogus isogeny theorem \cite{BOgusNotes}, plays a critical organizational role in capturing the additional structures on $\Prism_{R/A}$ (such as the Nygaard filtration). 

The paper \cite{BLMDRW} shows that when $(A,I)$ is a perfect crystalline prism (e.g., $(\mathbf{Z}_p,(p))$), one can reconstruct the de Rham--Witt complex $W\Omega^\bullet_R$ from the pair $(\Prism_{R/A},\widetilde{\varphi}_{R/A})$ by a pure homological algebra construction dubbed ``saturation''. Moreover, this construction has the potential to offer a better behaved alternative to the de Rham complex for singular varieties in characteristic $p$, analogous to the du Bois complex in characteristic $0$; we refer to \cite{IllusiedRWnew,OgusSaturateddRW} for more on these developments.
\end{remark}

\begin{remark}[Logarithmic analogs]
	The smoothness assumption on $X$ in Theorem~\ref{PrismaticCoh} is a ``good reduction'' hypothesis. While adequate for several purposes, this is often too restrictive for studying the generic fibre:  not every proper smooth scheme $X_\eta/C$ admits such a smooth model $X/\mathcal{O}_C$. A more natural assumption --- one that is conjecturally always satisfied, up to replacing models --- would be a form of logarithmic smoothness of $X/\mathcal{O}_C$ (e.g., semistability) in the sense of log geometry \cite{KatoLog}. Thus, one wants a version of \cite{BMS1,BhattScholzePrisms} in the logarithmic setting. This has been accomplished in \cite{CesKoshAinf,KoshLogPrism1,KoshYaoLogPrism2}; it is also possible to approach this problem by reduction to the smooth case using the language of infinite root stacks, following ideas of Olsson \cite{OlssonLogarithmic} (work in progress with Mathew).
\end{remark}

\begin{remark}[Non-abelian $p$-adic Hodge theory]
\label{rmk:NApHT}
Fix a bounded prism $(A,I)$ and a smooth $p$-adic formal $\overline{A}$-scheme $X$. Motivated by the precise form of Theorem~\ref{PrismaticCoh}, define a {\em prismatic $F$-crystal} on $(X/A)$ to be a vector bundle $\mathcal{E}$ on $((X/A)_\Prism, \mathcal{O}_\Prism)$ equipped with a Frobenius structure $\varphi_{\mathcal{E}}:\varphi^* \mathcal{E}[\frac{1}{I_\Prism}] \simeq \mathcal{E}[\frac{1}{I_\Prism}]$; see Definition~\ref{def:Crystals} for a more explicit description in a variant context. Prismatic $F$-crystals provide a viable notion of ``coefficients'' in the theory, somewhat analogous to the role played by harmonic bundles in complex non-abelian Hodge theory \cite{SimpsonICM,SimpsonHiggs}. In particular, given such an $F$-crystal $(\mathcal{E},\varphi_{\mathcal{E}})$, the specialization functors used in Theorem~\ref{PrismaticCoh} yield a $\mathbf{Z}_p$-local system $T(\mathcal{E})$ on the rigid generic fibre $X_\eta$ when $A$ is perfect, a vector bundle $\mathcal{E}_{\mathrm{dR}}$ with flat connection on $X/\overline{A}$,  an $F$-crystal $\mathcal{E}_{\mathrm{crys}}$ on $X \otimes_{\overline{A}} (\overline{A}/p)_{\mathrm{perf}}$, and (under certain auxiliary lifting data) a Higgs bundle $\mathcal{E}_{\mathrm{Higgs}}$ on $X/\overline{A}$. The relationship realized by these functors is rather close and has been investigated by various authors (such as \cite{MorrowTsuji,TianFCrys,BhattLurieAPC}). When $(A,I)$ is a Fontaine prism, this relationship is part of the $p$-adic Simpson correspondence pioneered by Faltings \cite{FaltingsSimpson1,FaltingsSimpson2,AbbesGrosTsujiSimpson}. On the other hand, if $(A,I)$ is a perfect crystalline prism, this relationship yields an alternative perspective on (at least the local aspects of) the non-abelian Hodge theory of \cite{OgusVologodsky}.
\end{remark}

\begin{remark}[Extension to the singular case via animation]
\label{rmk:Animation}
For several applications including most results discussed in this paper, it is important to extend the prismatic cohomology construction $R \mapsto \Prism_{R/A}$ (see Remark~\ref{rmk:DRW}) to possibly singular $\overline{A}$-algebras $R$. Directly imitating Definition~\ref{def:RelPrismSite} does not produce a computable or useable result. Instead, inspired by the construction of the cotangent complex and derived de Rham cohomology \cite{IllusieCC1,IllusieCC2} as well as their utility in a wide variety of problems \cite{Beilinsonpadic,BhattCompletions,Bhattpadic,BBLSZStabilization,IllusiePartial,GuoLiPeriod}, one extends the functor $\Prism_{-/A}$ to arbitrary $p$-complete $\overline{A}$-algebras by Quillen's non-abelian derived functor machinery \cite{QuillenHA} (dubbed {\em animation} by Clausen \cite{CesnaviciusScholzePurity}) as reformulated in \cite{LurieHTT}. The resulting complex $\Prism_{R/A}$ can be fairly efficiently controlled using the cotangent complex $L_{R/\overline{A}}$ thanks to the animated Hodge--Tate comparison, which makes this extension quite useable.
\end{remark}

\section{Absolute prismatic cohomology}
\label{sec:APC}

In \S \ref{sec:RelativePrismatic}, we fixed a base prism $(A,I)$ and discussed results about the relative prismatic cohomology of a smooth $p$-adic formal $\overline{A}$-scheme $X$. In this section, we describe the picture that arises if one does not fix a base prism $(A,I)$. This distinction is analogous to that between geometric and absolute \'etale cohomology in arithmetic, or that between singular cohomology and Deligne-Beilinson cohomology in Hodge theory. The objects considered here are newer than those in \S \ref{sec:RelativePrismatic}; consequently, some results are surely not optimal, and we have tried to indicate some natural further directions in the exposition.

\subsection{Definition and key examples}

To begin, let us recall the definition of the absolute prismatic site (obtained roughly from Definition~\ref{def:RelPrismSite} by discarding $(A,I)$).

\begin{definition}[The absolute prismatic site]
\label{Def:APS}
Given a $p$-adic formal scheme $X$, its absolute prismatic site $X_\Prism$ is the category of all bounded prisms $(B,J)$ equipped with a map $\mathrm{Spf}(B/J) \to X$, topologized using the flat topology; write $\mathcal{O}_\Prism$, $I_\Prism$ and $\overline{\mathcal{O}}_\Prism$ for the sheaves obtained by remembering $B$, $J$ or $B/J$ respectively. Write $R\Gamma_\Prism(X) := R\Gamma(X_\Prism,\mathcal{O}_\Prism) \ \rotatebox[origin=c]{-270}{$\circlearrowright$} \ \varphi$ and $R\Gamma_{\overline{\Prism}}(X) := R\Gamma(X_\Prism, \overline{\mathcal{O}}_\Prism)$ for the resulting cohomology theories. 
\end{definition}

If there exists a perfect prism $(A,I)$ and a map $X \to \mathrm{Spf}(\overline{A})$, the natural map $(X/A)_\Prism \to X_\Prism$ is an equivalence, so Theorem~\ref{PrismaticCoh} describes $R\Gamma_\Prism(X)$  in this case, e.g., $R\Gamma_\Prism(\mathrm{Spf}(\overline{A})) \simeq A$.  At the the other end, $\mathrm{Spf}(\mathbf{Z}_p)_\Prism$ is the opposite of the category of all bounded prisms. As this category has no final object, the cohomology theory $R\Gamma_\Prism(\mathrm{Spf}(\mathbf{Z}_p))$ is potentially interesting; in fact, we shall see in \S \ref{sec:K} that $R\Gamma_\Prism(\mathrm{Spf}(\mathbf{Z}_p))$  is closely related to the $p$-completed  algebraic $K$-theory of $\mathbf{Z}_p$. 

In this section, we shall be interested in the following objects on $X_\Prism$:

\begin{definition}[Crystals]
\label{def:Crystals}
Fix a $p$-adic formal scheme $X$.  A {\em prismatic crystal} (resp. {\em Hodge--Tate crystal}) $\mathcal{E}$ of vector bundles on $X$ is given by an assignment 
\begin{align*}
(B,J) \in X_\Prism \mapsto& \mathcal{E}(B) \in \mathrm{Vect}_B := \{ \text{finite projective } B\text{-modules}\} \\
 \Big(\text{resp. } (B,J) \in X_\Prism \mapsto& \mathcal{E}(B) \in \mathrm{Vect}_{B/J}\Big)
\end{align*}
that is compatible with base change in $(B,J) \in X_\Prism$.  A {\em prismatic $F$-crystal} of vector bundles on $X$ is given by a prismatic crystal $\mathcal{E}$ with an isomorphism $\varphi_{\mathcal{E}}:\varphi^* \mathcal{E}[\frac{1}{I_\Prism}] \simeq \mathcal{E}[\frac{1}{I_\Prism}]$ of $\mathcal{O}_\Prism$-modules. Similarly, one has analogous notions of crystals of perfect (or just $(p,I_\Prism)$-complete) complexes. 
\end{definition}

As in the relative case (Remark~\ref{rmk:NApHT}), there are realization functors carrying a prismatic $F$-crystal $\mathcal{E}$ on $X$ to a $\mathbf{Z}_p$-local system $T(\mathcal{E})$ on the rigid generic fibre $X_\eta$, a vector bundle $\mathcal{E}_{\mathrm{dR}}$ with flat connection on $X$, and an $F$-crystal $\mathcal{E}_{\mathrm{crys}}$ on $X \otimes_{\mathbf{Z}_p} \mathbf{F}_p$. The simplest examples of such crystals are as follows:

\begin{example}[Breuil-Kisin twists]
\label{ex:BKtwist}
For any prism $(B,J)$, one has a naturally defined invertible $B$-module $B\{1\}$ given heuristically by
\[ B\{1\} := J \otimes \varphi^* J \otimes \varphi^{2,*} J \otimes ....\]
This $B$-module comes equipped with a natural isomorphism $\varphi_{B\{1\}}:\varphi^* B\{1\}\simeq J^{-1} \otimes B\{1\}$, so the assignment $(B,J) \mapsto (B\{1\},\varphi_{B\{1\}})$ gives a prismatic $F$-crystal $(\mathcal{O}_\Prism\{1\}, \varphi_{\mathcal{O}_\Prism}\{1\})$ on $\mathrm{Spf}(\mathbf{Z}_p)_\Prism$ (and thus on $X_\Prism$ for any $X$); we refer to this $F$-crystal as the (first) Breuil-Kisin twist. The \'etale realization of $\mathcal{O}_\Prism\{1\}$ is identified with the usual Tate twist $\mathbf{Z}_p(1)$.
\end{example}

\begin{example}[Gauss-Manin $F$-crystals]
\label{ex:GM}
Fix a proper smooth map $f:Y \to X$ of $p$-adic formal schemes. The formalism of relative prismatic cohomology yields an $F$-crystal $Rf_* \mathcal{O}_\Prism$ of perfect complexes on $X_\Prism$: its value on a prism $(B,J) \in X_\Prism$ identifies with the relative prismatic complex $R\Gamma_\Prism( (Y \times_X \mathrm{Spf}(B/J))/B)$. Similarly, one obtains a Hodge--Tate crystal $Rf_* \overline{\mathcal{O}}_\Prism$ of perfect complexes on $X_\Prism$. The formation of $Rf_* \mathcal{O}_\Prism$ (resp. $Rf_* \overline{\mathcal{O}}_\Prism$) is compatible with the aforementioned realization functors.  Moreover, if $Y=\mathbf{P}^1 \times X$, then the prismatic logarithm  \cite{BhattLurieAPC} yields a natural isomorphism $\mathcal{H}^2(Rf_* \mathcal{O}_\Prism) \simeq \mathcal{O}_\Prism\{-1\}$ of $F$-crystals, giving a geometric description of the Breuil-Kisin twist.
\end{example}

\subsection{Hodge--Tate crystals}

In this subsection, we fix a perfect field $k$ of characteristic $p$, and write $W(k)_\Prism = \mathrm{Spf}(W(k))_\Prism$ for the absolute prismatic site of $W(k)$. Our goal is to explicitly describe the structure of Hodge--Tate crystals on $W(k)_\Prism$; we then  specialize this description to the Gauss-Manin case to obtain geometric consequences. For the former, we have:

\begin{proposition}[Sen theory, \cite{DrinfeldPrismatization,Drinfeld1dim,BhattLurieAPC}]
\label{prop:Sen}
The $\infty$-category $\widehat{\mathcal{D}}_{\mathrm{crys}}(W(k)_\Prism, \overline{\mathcal{O}}_\Prism)$ of Hodge--Tate crystals $\mathcal{E}$ of $p$-complete complexes on $W(k)_\Prism$ can be identified as the $\infty$-category of pairs $(E,\Theta)$ consisting of a $p$-complete object $E \in \mathcal{D}(W(k))$ and an endomorphism $\Theta:E \to E$ such that $\Theta^p-\Theta$ is locally nilpotent on $H^*(E/p)$; we refer to such pairs $(E,\Theta)$ as {\em Sen complexes} and $\Theta$ as the {\em Sen operator}.
\end{proposition}

The implicit functor carrying the crystal $\mathcal{E}$ to $E \in \mathcal{D}(W(k))$ in Proposition~\ref{prop:Sen} is given by evaluating at the object of $W(k)_\Prism$ obtained by base changing to $W(k)$ the $\mathbf{F}_p^*$-fixed points of the $q$-de Rham prism (Remark~\ref{rmk:qdR}).

\begin{remark}[The stacky approach to prismatic crystals, \cite{DrinfeldPrismatization,BhattLurieAPC}]
\label{rmk:CartWitt}
Proposition~\ref{prop:Sen} is proven via a stacky approach to prismatic cohomology, developed independently in \cite{DrinfeldPrismatization} (with a precursor in \cite{DrinfeldStacky}) and \cite{BhattLurieAPC}.  Using a tiny amount of derived algebraic geometry \cite{LurieSAG}, these works attach a stack $\mathrm{WCart}_X$ --- the {\em Cartier-Witt stack of $X$} (called the {\em prismatization $X^\Prism$} in \cite{DrinfeldPrismatization}) --- on $p$-nilpotent test rings to any $p$-adic formal scheme $X$. This stack comes equipped with an effective Cartier divisor $\mathrm{WCart}_X^{\mathrm{HT}} \subset \mathrm{WCart}_X$ called the {\em Hodge--Tate locus}. These stacks are devised to geometrize the study of crystals on the prismatic site: for quasi-syntomic $X$, there is a natural $\otimes$-identification of the $\infty$-category $\widehat{\mathcal{D}}_{\mathrm{crys}}(X_\Prism, \mathcal{O}_\Prism)$ of crystals of $(p,I_\Prism)$-complete complexes on $(X_\Prism,\mathcal{O}_\Prism)$ with the quasi-coherent derived $\infty$-category $\mathcal{D}_{qc}(\mathrm{WCart}_X)$; similarly the $\infty$-category $\widehat{\mathcal{D}}_{\mathrm{crys}}(X_\Prism, \overline{\mathcal{O}}_\Prism)$ of crystals of $p$-complete complexes on $(X_\Prism,\overline{\mathcal{O}}_\Prism)$ identifies with the quasi-coherent derived $\infty$-category $\mathcal{D}_{qc}(\mathrm{WCart}_X^{\mathrm{HT}})$. Proposition~\ref{prop:Sen} is then deduced from an explicit description of $\mathrm{WCart}_{W(k)}^{\mathrm{HT}}$ as $BG$ for a group scheme $G/W(k)$ whose representations are identified with Sen complexes.
\end{remark}

\begin{notation}[Diffracted Hodge cohomology]
\label{not:DiffHod}
Let $f:X \to \mathrm{Spf}(W(k))$ be a smooth map of $p$-adic formal schemes. Write $(R\Gamma(X,\Omega_X^{\DHod}),\Theta)$ for the Sen complex corresponding to the Hodge--Tate crystal $Rf_* \overline{\mathcal{O}}_\Prism \in \widehat{\mathcal{D}}_{\mathrm{crys}}(W(k)_\Prism, \overline{\mathcal{O}}_\Prism)$ via Proposition~\ref{prop:Sen}; we call $R\Gamma(X,\Omega_X^{\DHod})$ the {\em diffracted Hodge complex} of $X$. 
\end{notation}

The next result says that $R\Gamma(X,\Omega_X^{\DHod})$ is a slightly twisted form of the Hodge cohomology complex $\bigoplus_i R\Gamma(X, \Omega^i_X[-i])$, justifying the name ``diffracted Hodge cohomology".

\begin{theorem}[The Sen structure of $\Omega_X^{\DHod}$, \cite{BhattLurieAPC}]
\label{thm:DiffHod}
Let $X/W(k)$ be a smooth $p$-adic formal scheme. Then the Sen complex $(R\Gamma(X,\Omega_X^{\DHod}),\Theta)$  has a natural multiplicative increasing {\em conjugate} filtration $\mathrm{Fil}_{\mathrm{conj}}^\bullet$ equipped with natural isomorphisms 
\[ \mathrm{gr}^i_{\mathrm{conj}} (R\Gamma(X,\Omega_X^{\DHod}), \Theta) \simeq (R\Gamma(X,\Omega^i_{X/W(k)})[-i], \Theta=-i)\] 
for all $i$. 
\end{theorem}

Theorem~\ref{thm:DiffHod} also shows that the assignment $U \mapsto R\Gamma(U,\Omega_{U}^{\DHod})$ patches to a perfect complex $\Omega_X^{\DHod}$ on $X$, justifying the notation $R\Gamma(X,\Omega_X^{\DHod})$.

\begin{remark}[Relation to classical Sen theory, \cite{BhattLurieAPC}]
Fix a proper smooth map $f:X \to \mathrm{Spf}(W(k))$ of $p$-adic formal schemes; write $K=W(k)[1/p]$, fix a completed algebraic closure $C/K$, and write $G_K$ for the absolute Galois group of $K$. Classical results in $p$-adic Hodge theory \cite{SenGaloisRep,FaltingspHT} show that for each $n \geq 0$, the $C$-semilinear $G_K$-representation $H^n(X_C,\mathbf{Q}_p) \otimes_{\mathbf{Q}_p} C$ comes equipped with a canonical semisimple endomorphism $\theta_n$ whose eigenvalue decomposition yields the Hodge--Tate decomposition: we have
\[ H^n(X_C,\mathbf{Q}_p) \otimes_{\mathbf{Q}_p} C \simeq \bigoplus_{i=0}^n H^{n-i}(X, \Omega^{i}_{X/W(k)}) \otimes_{W(k)} C(-i),\]
with $\theta_n$ acting by $-i$ on the $i$-th summand on the right. Using the comparison isomorphisms in Theorem~\ref{PrismaticCoh}, one can roughly regard Theorem~\ref{thm:DiffHod} as an integral lift of this assertion: the value of $Rf_* \overline{\mathcal{O}}_\Prism$ on the Fontaine prism for $\mathcal{O}_C$ recovers $R\Gamma(X_C,\mathbf{Q}_p) \otimes_{\mathbf{Q}_p} C$ on inverting $p$, the Sen operator $\Theta$ from Theorem~\ref{thm:DiffHod} induces the Sen operator $\theta_n$ on each $H^n$ with the conjugate filtration from Theorem~\ref{thm:DiffHod} yielding the Hodge--Tate decomposition. It was a pleasant surprise to the author that the Sen operator admits a nice integral form.
\end{remark}

\begin{remark}[Drinfeld's refinement of Deligne-Illusie, \cite{DrinfeldPrismatization,BhattLurieAPC}]
In the setup of Theorem~\ref{thm:DiffHod}, there is a natural identification $R\Gamma(X,\Omega_X^{\DHod})/p \simeq R\Gamma_{\mathrm{dR}}(X_k)$ compatible with the conjugate filtration via Theorem~\ref{PrismaticCoh}. Drinfeld observed that the Sen operator then yields interesting consequences for $R\Gamma_{\mathrm{dR}}(X_k)$. More precisely, there is a $\mathbf{Z}/p$-grading on $R\Gamma_{\mathrm{dR}}(X_k)$ corresponding to the generalized eigenspace decomposition for the Sen operator $\Theta$, and the $i$-th conjugate graded piece $\mathrm{gr}_i^{\mathrm{conj}} \simeq R\Gamma(X_k,\Omega^i_{X_k})[-i]$ contributes only to the generalized eigenspace for $\Theta=-i$ by Theorem~\ref{thm:DiffHod}. In particular, if $\dim(X_k) < p$,  the conjugate filtration on $R\Gamma_{\mathrm{dR}}(X_k)$ splits canonically. This gives a conceptual new proof --- in fact a refinement --- of the seminal Deligne-Illusie result \cite{DeligneIllusie} on Hodge-to-de Rham degeneration (itself inspired by \cite{KatoVanishing,FontaineMessing,FaltingspHT}). As in \cite{DeligneIllusie},  one only needs a $W_2(k)$-lift of $X_k$ to obtain the Sen operator --- and thus the $\mathbf{Z}/p$-grading --- on $R\Gamma_{\mathrm{dR}}(X_k)$; this follows from  an analysis of $\mathrm{WCart}_{W_2(k)}^{\mathrm{HT}}$ similar to Proposition~\ref{prop:Sen}.  The results discussed in this paragraph  refine those in \cite{AchingerSuh} by one cohomological degree; another stacky proof was recently found in \cite{LiMondalEnd}.

For $X/W(k)$ a smooth formal scheme without any constraints on $\dim(X_k)$, one now obtains a residual nilpotent operator $\Theta+i$ on the generalized $\Theta$-eigenspace $R\Gamma_{\mathrm{dR}}(X_k)_i  \subset R\Gamma_{\mathrm{dR}}(X_k)$ corresponding to the eigenvalue $-i$; this operator seems to be a new piece of structure that awaits further investigation. 
\end{remark}

\subsection{The Nygaard filtration}

The absolute prismatic cohomology\footnote{The complex $R\Gamma_\Prism(X)$ as defined in Definition~\ref{Def:APS} works well under mild assumptions on the singularities of $X$ (e.g., for lci $X$). In general, one modifies the definition of $R\Gamma_\Prism(X)$ by a categorical procedure involving quasi-syntomic descent and animation; we do not elaborate on this further in this survey and refer to \cite{AMMNBFS,BhattLurieAPC} for more.} $R\Gamma_\Prism(X)$ of a $p$-adic formal scheme $X$ carries an important filtration $\mathrm{Fil}^\bullet_N$, called {\em the Nygaard filtration}. This filtration plays roughly the same role for prismatic cohomology as the Hodge filtration does for de Rham cohomology.  Moreover, for applications to algebraic topology (such as Theorem~\ref{OddVan} below), it is critical to understand this filtration. Its defining feature is that the Frobenius $\varphi$ on $R\Gamma_\Prism(X)$ carries $\mathrm{Fil}^i_N R\Gamma_\Prism(X)$ to $R\Gamma(X_\Prism, I_\Prism^i)$ for all $i$.  The relative version of this filtration is well understood, at least on graded pieces, thanks to the isomorphism \eqref{eq:Leta} and the Beilinson $t$-structure (see \cite[Proposition 5.8]{BMS2}). For the absolute version, one has a similar description:

\begin{theorem}[The Nygaard fibre sequence, \cite{BhattLurieAPC}]
\label{thm:NygaardFib}
For any $p$-adic formal scheme $X$ and any integer $i \geq 0$, there are natural fibre sequences
\begin{equation}
\label{eq:NygFib}
\tag{Nyg}
 \mathrm{gr}^i_N R\Gamma_\Prism(X) \to \mathrm{Fil}_i^{\mathrm{conj}} R\Gamma(X,\Omega_X^{\DHod}) \xrightarrow{\Theta+i} \mathrm{Fil}_{i-1}^{\mathrm{conj}} R\Gamma(X,\Omega_X^{\DHod}) 
 \end{equation}
 and
\begin{equation}
\tag{HT}
\label{eq:NygFib1}
R\Gamma_{\overline{\Prism}}(X)\{i\} := R\Gamma(X_\Prism, \overline{\mathcal{O}}_\Prism\{i\}) \to R\Gamma(X,\Omega_X^{\DHod}) \xrightarrow{\Theta+i} R\Gamma(X,\Omega_X^{\DHod}) 
\end{equation}
 with the convention that $\mathrm{Fil}_{<0}^{\mathrm{conj}} = 0$.
\end{theorem}

\begin{remark}[Calculations via the Nygaard fibre sequence]
\label{rmk:NygSeqCalc}
The sequence \eqref{eq:NygFib} is quite useful for calculations of the Nygaard filtration. For instance, in conjunction with the $\mathrm{THH}(-)$ variant of Theorem~\ref{KMot} below (see Remark~\ref{rmk:THHBMS2Variant} as well), one may use \eqref{eq:NygFib} to calculate $\pi_* \mathrm{THH}(R;\mathbf{Z}_p)$ for a $p$-completely smooth $\mathcal{O}_K$-algebra $R$, where $K$ is a discretely valued extension of $\mathbf{Q}_p$ with perfect residue field; this recovers calculations of \cite{BokstedtTHHZ,LindenstraussMadsen}. Comparing \eqref{eq:NygFib} and \eqref{eq:NygFib1} also quantifies the failure of the Frobenius map $\varphi: \mathrm{gr}^i_N R\Gamma_\Prism(X)  \to R\Gamma_{\overline{\Prism}}(X)\{i\}$ to be an isomorphism in terms of coherent cohomology, thus giving a new mechanism to study the so-called Segal conjecture for $\mathrm{THH}$.
\end{remark}

\begin{remark}[View $\mathrm{Spec}(\mathbf{Z})$ as a curve]
Several results in mathematics have been inspired by the seemingly nonsensical idea that $\mathrm{Spec}(\mathbf{Z})$ is a curve over some non-existent base $\mathbf{F}$. In $p$-adic arithmetic geometry, this idea can sometimes lead to useful (and testable!) predictions in conjunction with the following related heuristics: 
\begin{itemize}
\item Perfectoid rings (e.g., finite fields) are formally \'etale over $\mathbf{F}$.
\item Topologically finite type regular $p$-complete $\mathbf{Z}_p$-algebras $R$ are smooth over $\mathbf{F}$ of relative dimension $\dim(R)$ (the Krull dimension).
\end{itemize}

We briefly discuss some examples of such predictions. First, if $R$ is a perfectoid ring, then the $p$-completion of $L_{R/\mathbf{Z}_p}$ identifies with $R[1]$, which is consistent with the prediction of a transitivity triangle of cotangent complexes for $\mathbf{F} \to \mathbf{Z} \to R$ and the heuristics above; this was already essentially observed in \cite{FaltingsKodaira}.

Next, Theorem~\ref{thm:NygaardFib} was partially conceived based on these heuristics: the fibre sequence \eqref{eq:NygFib} is obtained as the associated graded of a fibre sequence of filtered complexes allowing one to compute the absolute Nygaard filtration in terms of the relative one; the underlying fibre sequence of complexes for the latter was guessed based on the analogy between $\mathbf{Z}$ and a smooth curve over a perfectoid ring.

Next, if one views the Hodge--Tate locus $\mathrm{WCart}_X^{\mathrm{HT}} \subset \mathrm{WCart}_X$ of the Cartier-Witt stack (Remark~\ref{rmk:CartWitt}) as a version of the Hodge stack (i.e., the classifying stack of the tangent bundle) over $\mathbf{F}$, then the second heuristic above predicts that $\mathrm{WCart}_X^{\mathrm{HT}}$ is well behaved if $X$ is regular, e.g., the map $\pi:\mathrm{WCart}_X^{\mathrm{HT}} \to X$ should be a gerbe, and $R\pi_*$ must have coherent cohomological dimension at most $\dim(X)$; the first of these predictions is true, while the second is true at least in dimension $1$ (\cite{BhattLurieAPC}). Relatedly, there are some recently defined candidate notions of differential forms relative to $\mathbf{F}$ (\cite{SaitoFW,HochsterJeffriesJacobian,DKRZ}); it would be interesting to find a direct connection between the stack $\mathrm{WCart}_X^{\mathrm{HT}}$ and these objects.

Finally, let us remark that the philosophy discussed in this remark also featured in Scholze's report for the previous ICM \cite{ScholzeICMRio}, and has paid amazing dividends in geometrizing the local Langlands correspondence in recent years \cite{SWBerkeley,ScholzeDiamonds,FarguesScholze}. 
\end{remark}

\begin{remark}[$p$-adic Tate twists, \cite{BMS2,BhattScholzePrisms}]
\label{rmk:TateTwist}
An early observable extracted from absolute prismatic cohomology was a good notion of $p$-adic Tate twists $\mathbf{Z}_p(i)(-)$ in mixed characteristic: these are functors on $p$-adic formal schemes $X$ defined by a fibre sequence
\begin{equation}
\label{eq:TateDef}
\tag{Syn}
 \mathbf{Z}_p(i)(X) \to \mathrm{Fil}^i_N R\Gamma_\Prism(X)\{i\} \xrightarrow{\varphi-1} R\Gamma_\Prism(X)\{i\}
 \end{equation}
for all $i \geq 0$. These functors are often called {\em syntomic complexes} for mixed characteristic rings as they extend those in characteristic $p$ considered in \cite{MilneDuality,KatoSyntomic}. One can also regard $\mathbf{Z}_p(i)(-)$ as as a form of \'etale motivic cohomology in weight $i$ (see the forthcoming Remark~\ref{rmk:EtaleMot}).  In fact, for formally smooth $\mathcal{O}_K$-schemes with $K/\mathbf{Q}_p$ finite, the syntomic complexes $\mathbf{Z}_p(i)(-)$ agree with the $p$-adic \'etale Tate twists of Geisser-Sato-Schneider \cite{GeisserMot,SchneiderTate, SatoTateTwist} defined by glueing motivic complexes on the generic and special fibres \cite{BCMMot}. We refer the reader to \cite{BMS2,AMMNBFS,BhattScholzePrisms,BhattLurieAPC} for more on these syntomic complexes.
\end{remark}

\begin{remark}[$p$-adic Picard and Brauer groups via coherent cohomology]
The syntomic complex from Remark~\ref{rmk:TateTwist} in weight $1$ has the following relationship with $\mathbf{G}_m$ (\cite[Proposition 7.17]{BMS2}), as motivic intuition predicts: for any $p$-adic formal scheme $X$, we have
\[ \mathbf{Z}_p(1)(X) \simeq R\Gamma(X_{et},\mathbf{G}_m)^{\wedge}[-1],\]
where the completion is $p$-adic. Plugging this into the sequence \eqref{eq:TateDef} gives a fibre sequence 
\begin{equation}
\label{eq:Zp1}
\tag{Lef}
  R\Gamma(X_{et},\mathbf{G}_m)^{\wedge}[-1] \to \mathrm{Fil}^1_N R\Gamma_\Prism(X)\{1\} \xrightarrow{\varphi-1} R\Gamma_\Prism(X)\{1\}
 \end{equation}
that can be regarded as a weak $p$-adic analog of the Lefschetz $(1,1)$-theorem, e.g., it enables one to compute the $p$-completion of $\mathrm{Pic}(X)$ or $\mathrm{Br}(X)$ in terms of absolute prismatic cohomology, and thus ultimately via coherent cohomology. 

The idea described in the previous paragraph inspired the eventual proof of (a generalization of) Gabber's purity conjectures for Picard and Brauer groups in \cite{CesnaviciusScholzePurity}. In a global direction, Cotner and Zavyalov have recently used \eqref{eq:Zp1} to prove the vanishing of $\mathrm{Pic}^{\tau}(X)$ for complete intersection surfaces $X \subset \mathbf{P}^N$ in characteristic $p$ (in progress), settling a question left open since \cite{SGA2}. In a different direction, the sequence \eqref{eq:Zp1} can be used to prove that $R \mapsto R\Gamma(\mathrm{Spf}(R)_{et}, \mathbf{G}_m)^{\wedge}$ commutes with sifted colimits in $R$ (in the $p$-complete world); this allows one to reduce general questions about $R\Gamma(\mathrm{Spf}(R)_{et}, \mathbf{G}_m)^{\wedge}$ to particularly nice rings, and played an important role in Bragg and Olsson's work \cite{BraggOlssonRep} on finiteness results for higher direct images of finite flat group schemes along projective morphisms in characteristic $p$.
 \end{remark}

\subsection{Galois representations}

In this subsection, fix a discretely valued field $K/\mathbf{Q}_p$ with perfect residue field $k$. We discuss the relationship of prismatic $F$-crystals over $X=\mathrm{Spf}(\mathcal{O}_K)$ and Galois representations of $G_K = \mathrm{Gal}(\overline{K}/K)$. 

For a prime $\ell \neq p$,  the notion of {\em unramifiedness} for $\mathbf{Z}_\ell$- or $\mathbf{Q}_\ell$-representations of $G_K$ is a Galois-theoretic analog of the property of having ``good reduction'' for varieties over $K$: viewed as an $\ell$-adic local system on $\mathrm{Spec}(K)$, an unramified $G_K$-representation is exactly one that extends to a local system over $\mathrm{Spf}(\mathcal{O}_K)$. In contrast, for $\mathbf{Z}_p$- or $\mathbf{Q}_p$-representations, unramifiedness is too restrictive: even the cyclotomic character --- or any nonzero $H^i(Y_{\overline{K}}, \mathbf{Q}_p)$ with $Y/\mathcal{O}_K$ smooth projective and $i > 0$ --- is not unramified. To remedy this, Fontaine invented \cite{FontaineAnnals} the notion of {\em crystalline} $G_K$-representations;  it has been stunningly successful at capturing the desired ``good reduction'' intuition. On the other hand, any prismatic $F$-crystal $\mathcal{E}$ on $\mathrm{Spf}(\mathcal{O}_K)$ gives rise to a $G_K$-representation $T(\mathcal{E})$ as well an $F$-crystal $\mathcal{E}_{\mathrm{crys}}$ on $k$ (see Definition~\ref{def:Crystals} and following discussion); these have the same rank, so one may view $\mathcal{E}_{\mathrm{crys}}$ as ``a special fibre'' of $T(\mathcal{E})$, suggesting that the prismatic $F$-crystal $\mathcal{E}$ itself should be viewed as a witness for a ``good reduction'' of $T(\mathcal{E})$. The following theorem shows that these two perspectives on good reduction for $p$-adic representations coincide:

\begin{theorem}[Prismatic $F$-crystals and crystalline $G_K$-representations, \cite{BhattScholzeFCrys}]
\label{thm:FCrys}
The \'etale realization functor $\mathcal{E} \mapsto T(\mathcal{E})$ gives an equivalence of the category  of prismatic $F$-crystals on $\mathrm{Spf}(\mathcal{O}_K)$ with the category of $\mathbf{Z}_p$-lattices in crystalline $\mathbf{Q}_p$-representations of $G_K$. 
\end{theorem}

Thus, prismatic $F$-crystals on $\mathrm{Spf}(\mathcal{O}_K)$ provide a reasonable notion for ``local systems on $\mathrm{Spf}(\mathcal{O}_K)$ with $\mathbf{Z}_p$-coefficients''.

\begin{remark}
Theorem~\ref{thm:FCrys} can be viewed as a refinement of Kisin's classification of crystalline $G_K$-representations \cite{Kisin}; in particular, this refinement attaches prismatic meaning to the integrality properties of a somewhat mysterious connection in \cite{Kisin}. An alternative proof of Theorem~\ref{thm:FCrys} was since given in \cite{DuLiuFCrys}, relying on the theory in \cite{LiuNoteLattice}; see also \cite{WuGalois}.
\end{remark}

Various results in the deformation theory of $G_K$-representations (e.g., \cite{KisinPst,EmertonGeeStack}) indicate it would be fruitful to extend the notion of crystalline $G_K$-representations to torsion coefficients or even to the derived category. However, as the property of  being crystalline is essentially a rational concept, it is not clear how to proceed. Theorem~\ref{thm:FCrys} points to a way forward, e.g., perhaps prismatic $F$-crystals with $\mathcal{O}_\Prism/p^n$-coefficients are a reasonable candidate for crystalline $\mathbf{Z}/p^n$-representations? While satisfactory for describing $\mathbf{Z}_p$-local systems, this approach does not quite lead to a reasonable derived theory as the definition of a prismatic $F$-crystal $(\mathcal{E},\varphi_{\mathcal{E}})$ is not quantitative enough: the isomorphism $\varphi_{\mathcal{E}}$ does not come equipped with bounds on its poles/zeroes, leading to certain poorly behaved $\mathrm{Ext}$-groups. Instead, the correct objects seem to be perfect complexes on an enlargement of the Cartier-Witt stack $\mathrm{WCart}_{\mathcal{O}_K}$ (Remark~\ref{rmk:CartWitt}) constructed by Drinfeld \cite{DrinfeldPrismatization}; we describe one piece of evidence for this correctness assertion in the rest of the subsection.

Write $\mathcal{D}^{\varphi}_{\mathrm{perf}}(\mathrm{WCart}^+_{\mathcal{O}_K})$ for the $\infty$-category of perfect complexes on the stack $\mathrm{Spf}(\mathcal{O}_K)^{\Prism ''}$ from \cite[\S 1.8]{DrinfeldPrismatization}; let us call such objects {\em prismatic $F$-gauges on $\mathcal{O}_K$}\footnote{The definition of $\mathrm{WCart}^+_{\mathcal{O}_K}$ in \cite{DrinfeldPrismatization} (denoted $\Sigma'_{\mathcal{O}_K}$ there) is inspired by the Fontaine--Jansen theory \cite{FontaineJannsen} of $F$-gauges in crystalline cohomology.}. Given such an $F$-gauge $\mathcal{E}$, write $R\Gamma^{\varphi}(\mathrm{WCart}_{\mathcal{O}_K}^+, \mathcal{E})$ for its global sections. To a first approximation, a prismatic $F$-gauge $\mathcal{E}$ consists of a prismatic $F$-crystal $E$ of perfect complexes on $\mathrm{WCart}_{\mathcal{O}_K}$ equipped with the additional datum of a Nygaard-style filtration on $R\Gamma(\mathrm{WCart}_{\mathcal{O}_K},E)$; in fact, this can be made precise if $\mathcal{O}_K$ is replaced by a qrsp ring (work in progress with Lurie). The prismatic $F$-crystals from Examples~\ref{ex:BKtwist} and \ref{ex:GM} have natural lifts to prismatic $F$-gauges. The promised piece of evidence is the following result:

\begin{theorem}[A Lagrangian property, \cite{BhattLurieAPG}]
\label{thm:Duality}
Assume $K$ is unramified. Let $\mathcal{E} \in \mathcal{D}^{\varphi}_{\mathrm{perf}}(\mathrm{WCart}^+_{\mathcal{O}_K})$ be a prismatic $F$-gauge with $\mathcal{O}$-linear dual $\mathbf{D}(\mathcal{E})$ and \'etale realization $T(\mathcal{E})$ in an appropriate derived category  of continuous $\mathbf{Z}_p$-representations of $G_K$. Then the natural map
\[ R\Gamma^{\varphi}(\mathrm{WCart}^+_{\mathcal{O}_K}, \mathcal{E}) \to R\Gamma(G_{K}, T(\mathcal{E})) \]
is the exact annihilator of the corresponding map of local Tate duals, i.e., there is a natural fibre sequence
\begin{equation}
\label{eq:Lag}
\tag{Lag}
 R\Gamma^{\varphi}(\mathrm{WCart}^+_{\mathcal{O}_K}, \mathcal{E}) \to R\Gamma(G_{K}, T(\mathcal{E})) \to \Big(R\Gamma^{\varphi}(\mathrm{WCart}^+_{\mathcal{O}_K}, \mathbf{D}(\mathcal{E})\{1\}[2])\Big)^\vee
 \end{equation}
where $(-)^\vee = \mathrm{RHom}_{\mathbf{Z}_p}(-,\mathbf{Z}_p)$ on the rightmost term.
\end{theorem}

Theorem~\ref{thm:Duality} is work in progress with Lurie \cite{BhattLurieAPG}; the statement is likely not quite optimal yet (e.g., we hope to show it for ramified $K$ as well). 

\begin{remark}[The crystalline part of Galois cohomology]
Given a $\mathbf{Q}_p$-representation $V$ of $G_K$,  Bloch--Kato constructed \cite[\S 3]{BlochKato} the ``crystalline part'' $H^1_f(G_K, V) \subset H^1(G_K, V)$ of the Galois cohomology of $V$, and proved that the crystalline parts for $V$ and $V^\vee(1)$ are orthogonal complements under the Tate duality pairing. When $V$ is crystalline, Theorem~\ref{thm:Duality} may be viewed as an integral and cochain-level variant of this statement. Such integral refinements have been formulated previously in special cases (e.g., \cite{GazakiFiner}); it would be interesting to compare them to Theorem~\ref{thm:Duality}.
\end{remark}

\begin{remark}[The analogy with $3$-manifolds]
In the Mazur(-Mumford) analogy between number rings and $3$-manifolds \cite{MazurAlexander,LLMazurPost}, the scheme $\mathrm{Spec}(K)$ corresponds to a Riemann surface $\Sigma$ while $\mathrm{Spf}(\mathcal{O}_K)$ corresponds to a $3$-manifold $M$ with boundary $\partial M = \Sigma$. A standard topological result states that the space $\mathrm{Loc}(\Sigma)$ of local systems on $\Sigma$ has a symplectic structure induced by Poincar\'e duality on $\Sigma$, and the restriction map $\mathrm{Loc}(M) \to \mathrm{Loc}(\Sigma)$ is Lagrangian (see \cite[Proposition 3.27]{FreedChernSimons}). The sequence \eqref{eq:Lag} may be viewed as an arithmetic analogue of the infinitesimal form of this result, with the role of a local system on $\mathrm{Spf}(\mathcal{O}_K)$ played by prismatic $F$-gauges; in fact, this picture motivated the discovery of Theorem~\ref{thm:Duality}.
\end{remark}

\begin{remark}[Lichtenbaum-Quillen for $\mathcal{O}_K$]
\label{rmk:LichQuillenOK}
The Breuil-Kisin prismatic $F$-gauges $\mathcal{O}\{i\}$ compute the $p$-adic Tate twists from Remark~\ref{rmk:TateTwist}, i.e, we have natural identifications $R\Gamma^{\varphi}(\mathrm{WCart}^+_{\mathcal{O}_K}, \mathcal{O}\{i\}) \simeq \mathbf{Z}_p(i)(\mathcal{O}_K)$. Using the vanishing of $\mathbf{Z}_p(i)(-)$ for $i < 0$, the sequence \eqref{eq:Lag}  then implies that the natural map
\[ \mathbf{Z}_p(i)(\mathcal{O}_K) \to R\Gamma(G_K, \mathbf{Z}_p(i))\]
is an equivalence for $i \geq 2$. Under the relationship of either side to the \'etale $K$-theory of $\mathcal{O}_K$ and $K$ as well as the localization sequence in $K$-theory, this result was essentially known (\cite{HesselholtMadsen}); nevertheless, Theorem~\ref{thm:Duality} provides a different conceptual explanation. 
\end{remark}

\section{Algebraic $K$-theory}
\label{sec:K}

Quillen's algebraic $K$-theory \cite{QuillenK1} functor attaches a space (in fact, a spectrum) $K(X)$ to a scheme $X$, generalizing the construction of the Grothendieck group $K_0(X)$; the study of these invariants and their generalizations is an important pursuit in modern algebraic topology. In fact, its impact extends far beyond algebraic topology: the higher $K$-groups $K_i(X)$ feature prominently in some of the deepest conjectures in arithmetic geometry. In this section, we report on some recent progress in understanding the structural features of the $p$-completed algebraic $K$-theory spectrum $K(R;\mathbf{Z}_p)$ of a $p$-complete ring $R$, with an emphasis on connections to prismatic cohomology; the case of $\ell$-adic completions for $\ell \neq p$ is classical, going back to work of Thomason \cite{ThomasonKEtale}, Suslin \cite{SuslinRigidK} and Gabber \cite{GabberRigidK}. More complete recent surveys of material covered in this section include \cite{HesselholtNikolausSurvey,MFOReportTC,MathewSurveyTHH}.

In classical algebraic topology, combining the Atiyah-Hirzeburch spectral sequence with Bott periodicity gives the following structure on the complex $K$-theory $K^{\mathrm{top}}(X)$ of a (reasonable) topological space $X$:

\[
\tag{$\mathrm{Filt}_K$}
\label{Ktop}
  \parbox{4.5in}{ The $K$-theory spectrum $K^{\mathrm{top}}(X)$ admits a natural filtration with $\mathrm{gr}^i$ identified with the shifted singular cohomology complex $R\Gamma(X,\mathbf{Z})[2i]$. \\ } 
  \]

 In recent work, motivated in part by conjectures of Beilinson \cite{BeilinsonHeight} and Hesselholt, a $p$-adic analog of \eqref{Ktop} for (\'etale sheafified) algebraic $K$-theory of $p$-complete rings has been established, with the role of singular cohomology now played by prismatic cohomology. To explain this better, recall that algebraic topologists often study algebraic $K$-theory through a ``cyclotomic trace'' map 
 \[ \mathrm{Tr}:K(-) \to \mathrm{TC}(-),\] 
 where $\mathrm{TC}(-)$ is the {\em topological cyclic homology} functor; this invariant of rings was invented by contemplating Hochschild homology over the sphere spectrum, goes back to \cite{BHM}, and was recently given a simple $\infty$-categorical definition in \cite{NikolausScholze}. The trace map is a powerful calculational tool (see \cite{HesselholtMadsenPolytope,HesselholtMadsenFinite,HesselholtMadsen,GeisserHessdRW} for some successes), and there are two main reasons for this. First,
$\mathrm{TC}(-)$ is built, via a rather elaborate homotopical procedure, from objects of coherent cohomology (namely, differential forms) and is thus potentially more computable than $K$-theory. Secondly, the trace map turns out to yield a very good approximation of $K$-theory in various situations; in our $p$-adic context, the state of the art is the following:
 
 \begin{theorem}[$p$-adic \'etale $K$-theory is $\mathrm{TC}$,  \cite{CMM,ClausenMathewHyper}]
 \label{thm:KTC}
For $p$-complete rings $R$, the trace map $K(R;\mathbf{Z}_p) \to \mathrm{TC}(R;\mathbf{Z}_p)$ identifies the target with the $p$-completed \'etale sheafified $K$-theory of $R$. Moreover, the \'etale sheafification is not necessary in sufficiently large degrees if $R$ satisfies mild finiteness conditions.
 \end{theorem}
 
Via Theorem~\ref{thm:KTC}, the following result can be viewed as a $p$-adic analog of the Atiyah-Hirzebuch part of \eqref{Ktop}:

\begin{theorem}[The motivic filtration on \'etale $K$-theory, \cite{BMS2,AMMNBFS}]
\label{KMot}
As a functor on $p$-complete rings, there is a natural ``motivic'' filtration on $\mathrm{TC}(-; \mathbf{Z}_p)$ with $\mathrm{gr}^i_{\mathrm{mot}} \mathrm{TC} (-;\mathbf{Z}_p)$ naturally identified with the shifted syntomic complex $\mathbf{Z}_p(i)(-)[2i]$ from Remark~\ref{rmk:TateTwist}.
\end{theorem}

\begin{remark}[Variants for $\mathrm{THH}$ and cousins]
\label{rmk:THHBMS2Variant}
Let us briefly recall the  \cite{NikolausScholze} approach to calculating $\mathrm{TC}$. For a commutative ring $R$, the $p$-completed topological Hochschild homology spectrum $\mathrm{THH}(R;\mathbf{Z}_p)$ comes equipped with a natural action of the circle $S^1$ and a certain Frobenius map. One can then define auxiliary invariants $\mathrm{TC}^{-}(R;\mathbf{Z}_p) := \mathrm{THH}(R;\mathbf{Z}_p)^{hS^1}$ and $\mathrm{TP}(R;\mathbf{Z}_p) := \mathrm{THH}(R;\mathbf{Z}_p)^{tS^1}$ together with two natural maps $\mathrm{can},\varphi:\mathrm{TC}^{-}(R;\mathbf{Z}_p) \to \mathrm{TP}(R;\mathbf{Z}_p)$. The paper \cite{NikolausScholze} then proves there is a natural fibre sequence 
\begin{equation}
\label{eq:NS}
\tag{TC}
 \mathrm{TC}(R;\mathbf{Z}_p) \to \mathrm{TC}^{-}(R;\mathbf{Z}_p) \xrightarrow{\varphi-\mathrm{can}} \mathrm{TP}(R;\mathbf{Z}_p),
 \end{equation}
 thereby yielding a clean modern construction of $\mathrm{TC}(R;\mathbf{Z}_p)$. 

The construction of the motivic filtration on $\mathrm{TC}(-;\mathbf{Z}_p)$ in Theorem~\ref{KMot} is sufficiently flexible to ensure that similar ideas also yield compatible ``motivic filtrations'' on $\mathrm{TC}^{-}(R;\mathbf{Z}_p)$ and $\mathrm{TP}(R;\mathbf{Z}_p)$. In fact, \cite{BMS2} lifts the sequence \eqref{eq:NS} to a filtered fibre sequence that recovers the sequence \eqref{eq:TateDef} on associated graded pieces, up to Nygaard completions. In particular, one recovers (Nygaard completed) absolute prismatic cohomology as the associated graded of a filtration on $\mathrm{TP}$.

Note that while $\mathrm{THH}(-)$ and cousins are non-commutative invariants (i.e., can be defined for arbitrary stable $\infty$-categories), the construction of the motivic filtration crucially uses algebraic geometry; it is unclear if analogous filtrations exist even in slightly more general settings, e.g., for $\mathrm{TC}(\mathcal{C})$ for a symmetric monoidal stable $\infty$-category $\mathcal{C}$.
\end{remark}

\begin{remark}[Comparison with the Hodge filtration on classical Hochschild homology]
The topological Hochschild homology of a commutative ring $R$ is defined as $\mathrm{THH}(R) := \mathrm{HH}(R/\mathbf{S})$, i.e., it is the Hochschild homology relative to the sphere spectrum $\mathbf{S}$. From this optic, Theorem~\ref{KMot} and the variants in Remark~\ref{rmk:THHBMS2Variant} are analogs of known constructions in the Hochschild homology of ordinary rings. For instance, given a smooth algebra $R$ over a commutative ring $k$, there is a natural filtration on $\mathrm{HC}^{-}(R/k) := \mathrm{HH}(R/k)^{hS^1}$ with $\mathrm{gr}^i$ identified with the Hodge filtration level $\mathrm{Fil}^i_{\mathrm{Hodge}} \Omega^\bullet_{R/k}[2i]$ (see \cite[\S 5.1.2]{LodayCyclic} and \cite{WeibelHodge} for characteristic $0$, and \cite[\S 5.2]{BMS2} and \cite{AntieauHP} in general). Specializing to $k = \overline{A}$ for a Fontaine prism $(A,I)$, this allows one to recover crystalline and de Rham cohomology --- but not \'etale cohomology of the generic fibre --- as graded pieces of a natural filtration on classical Hochschild homology and its variants. Theorem~\ref{KMot} and the variants in Remark~\ref{rmk:THHBMS2Variant} thus contain the surprise that working relative to the sphere spectrum permits one to see the \'etale cohomology of the generic fibre as well: one can in fact recover prismatic cohomology.
\end{remark}

\begin{remark}[Origin story]
\label{rmk:OriginKMot}
Let $(A,I)$ be a Fontaine prism (Example~\ref{PerfPrism}). Write $A_{\mathrm{crys}}$ for the ring obtained by formally adjoining divided powers of $I$ to $A$ in the $p$-complete setting. Given a smooth proper scheme $X$ over $\mathcal{O}_C = \overline{A}$, its absolute crystalline cohomology $R\Gamma_{\mathrm{crys}}(X)$ is a perfect complex of $A_{\mathrm{crys}}$-modules with a Frobenius structure. The relation between the theory of Breuil modules \cite{BreuilMod} and Breuil-Kisin modules \cite{Kisin} in Galois representation theory strongly suggested that $R\Gamma_{\mathrm{crys}}(X)$ ought to descend naturally along $A \to A_{\mathrm{crys}}$. Separately, Hesselholt had calculated  \cite{HesselholtCp} that $\pi_0 \mathrm{TP}(\mathcal{O}_C;\mathbf{Z}_p)$ equals $A$. Comparing this to the known fact that $\pi_0 \mathrm{HP}(\mathcal{O}_C;\mathbf{Z}_p)$ equals $A_{\mathrm{crys}}$ (up to a completion), it was natural to speculate that for any $X/\mathcal{O}_C$ as above, one could find a filtration on $\mathrm{TP}(X;\mathbf{Z}_p)$ whose graded pieces realize the desired descent of $R\Gamma_{\mathrm{crys}}(X)$ along $A \to A_{\mathrm{crys}}$; this eventually led to the $\mathrm{TP}(-)$-variant of Theorem~\ref{KMot} and gave a  construction of prismatic cohomology over the Fontaine and Breuil-Kisin prisms \cite[\S 11]{BMS2}. In fact, as $\mathrm{TP}(-;\mathbf{Z}_p)$ is independent of the base $\mathcal{O}_C$, this also gave the first construction of absolute prismatic cohomology \cite[\S 7.3]{BMS2}.
\end{remark}

\begin{remark}[\'Etale motivic cohomology]
\label{rmk:EtaleMot}
We briefly explain why Theorem~\ref{KMot} can be viewed as a $p$-completed and \'etale sheafified analog of the filtration of algebraic $K$-theory by motivic cohomology (defined via Bloch's higher Chow groups \cite{BlochHigherChow}). Recall that the latter geometric motivic filtration was conjectured to exist in \cite{BeilinsonHeight}, and established in many cases, including most smooth cases, in \cite{BlochLichtenbaumSS,FriedSusK,VoevodskyOpenSlice,LevineCon}; see \cite{HoyoisHopkinsMorel} for a clean construction. It is thus natural to conjecture that the $p$-completed \'etale sheafification of this geometric motivic filtration identifies with the one in Theorem~\ref{KMot}.  For smooth varieties over a perfect field $k$ of characteristic $p$, this is indeed the case by combining \cite{GeisserLevine} and \cite{BMS2}. In mixed characteristic, while we do not know the full story, a positive answer at the associated graded level  follows from the comparison result from \cite{BCMMot} mentioned in Remark~\ref{rmk:TateTwist}. Let us also remark that \cite{MorrowLuders} has established the expected  relationship of Milnor $K$-theory (extended as in \cite{KerzMilnor}) and the $(i,i)$-part of the syntomic complexes, giving a $p$-adic analog of \cite{NesterenkoSuslin,TotaroMilnor}. 

The above discussion raises a natural question: as Theorem~\ref{KMot} applies to any $p$-complete ring, can the domain of definition of the geometric motivic filtration of the previous paragraph also be extended to all $p$-complete rings? In particular, is there a meaningful geometric motivic filtration on $K(R)$ for a non-reduced $p$-complete ring $R$? Thanks to a forthcoming result of Mathew on glueing the filtration from Theorem~\ref{KMot} with the \'etale sheafified Postnikov filtration on $K(1)$-local $K$-theory, there is a variant of Theorem~\ref{KMot} for any ring, so one may even reasonably ask these questions for all rings.
\end{remark}

\begin{remark}[Constructing the motivic filtration via quasisyntomic descent]
\label{rmk:QSynDescent}
The construction of the motivic filtration in Theorem~\ref{KMot} is quite different from that of the geometric motivic filtration mentioned in Remark~\ref{rmk:EtaleMot}. Indeed, the general case of Theorem~\ref{KMot} is proven in \cite{AMMNBFS} (see also \cite{BhattLurieAPC}) by reducing (via animation as in Remark~\ref{rmk:Animation}) to the quasisyntomic case treated in \cite{BMS2}. The latter has two essential ingredients. The topological ingredient is B\"{o}kstedt's fundamental periodicity result \cite{BokstedtTHH} that $\pi_* \mathrm{THH}(\mathbf{F}_p) = \mathbf{F}_p[u]$ for a degree $2$ class $u$; see \cite{KrauseNikolaus} for a quick modern proof based on properties of the dual Steenrod algebra,  \cite{KaledinFonarevBokstedt} for an overview of other approaches, and \cite{HesselholtNikolausSurvey} for a deduction of Bott periodicity from B\"{o}kstedt periodicity. The new algebraic input is the flat descent property (\cite{BhattCompletions,BMS2}) of the cotangent complex, used in conjunction with the very perfectoid idea (going back in spirit to \cite{FontaineMessing}) that working with certain infinitely ramified covers can  ``discretize''  constructions involving differential forms in the $p$-adic world.
\end{remark}

As Remark~\ref{rmk:OriginKMot} explains, the first construction of absolute prismatic cohomology was through Theorem~\ref{KMot} and variants. However, thanks to the alternative and more direct construction via the prismatic site, one can now use Theorem~\ref{KMot} as a tool to study $K$-theory via prismatic cohomology. For instance,  this approach gives the following result: 
  
\begin{theorem}[The odd vanishing theorem, \cite{BhattScholzePrisms}]
\label{OddVan}
For odd $i$, the functor $\pi_i K(-;\mathbf{Z}_p)$ is quasi-syntomic locally $0$ on the category of quasi-syntomic rings.
\end{theorem}

Theorem~\ref{OddVan} can be regarded as a variant of the Bott periodicity part of $\eqref{Ktop}$ in the algebraic setting: while periodicity is known to be false due to geometric phenomena, we still have vanishing in odd degrees. The proof in \cite{BhattScholzePrisms} relies on Andr\'e's flatness lemma \cite{AndreDSC}, and it would be interesting to find a more explicit description of the necessary covers. 

\begin{remark}[Further relations to $p$-adic arithmetic geometry]
Another application of prismatic cohomology to $K$-theory was  proving that $L_{K(1)} K(R) \simeq L_{K(1)} K(R[1/p])$ for any associative ring $R$ (\cite{BCMK1}). This equality is a $K$-theoretic avatar of the \'etale comparison from \cite{BhattScholzePrisms} and was proved in \cite{BCMK1} via explicit calculations in prismatic cohomology; it has since been reproved and significantly extended using purely homotopy-theoretic methods in \cite{LMMT}. 

In the reverse direction (and preceding most of the developments reported in this paper), \cite{BlochEsnaultKerz} used results from topological cyclic homology \cite{GeisserHesselholtRelative} to prove the infinitesimal portion of the $p$-adic variational Hodge conjecture in the unramified case. The extension to the ramified case was recently obtained in \cite{AMMNBFS} as a consequence of a purely $K$-theoretic assertion called the Beilinson fibre square. Using this square and Theorem~\ref{KMot}, \cite{AMMNBFS} also gave a simple description of the rationalized syntomic complexes $\mathbf{Z}_p(i)(-)[1/p]$ via derived de Rham  cohomology. This description is quite useful as derived de Rham cohomology is more computable in practice than prismatic cohomology; in fact, this description formed  an essential ingredient in the classification of crystalline representations given in Theorem~\ref{thm:FCrys}.

The connections discussed above have mostly concerned relative prismatic cohomology. It seems likely that a better understanding of absolute prismatic cohomology (as in \S \ref{sec:APC}) will lead to more refined applications. For instance, \cite{LiuWangTC} recovers rather conceptually the highly non-trivial calculation  \cite{HesselholtMadsen} of the $K$-theory of local fields $K/\mathbf{Q}_p$ by exploiting certain covers of the final object in the absolute prismatic topos of $\mathcal{O}_K$ coming from Breuil-Kisin prisms. Other related observations are discussed in Remark~\ref{rmk:NygSeqCalc} and Remark~\ref{rmk:LichQuillenOK}.

\end{remark}

\section{Commutative algebra and birational geometry}
\label{sec:CABG}

The Kodaira vanishing theorem (as well as the generalization by Kawamata--Viehweg) is one of the most important foundational results in complex algebraic geometry; it is especially useful in birational geometry. Its (original) proof relies crucially on Hodge theory, and thus no longer applies in positive/mixed characteristic. In fact, the result is known to be false in those settings \cite{RaynaudKodaira}; alongside the non-availability of resolution of singularities in dimensions $\geq 4$, this is a major obstacle  to progress in birational geometry in postive/mixed characteristic. About a decade ago, Schwede observed \cite{SchwedeCanonical} that methods from $F$-singularity theory in positive characteristic commutative algebra can sometimes be used as a substitute for the use of vanishing theorems in positive characteristic algebraic geometry;  this eventually led to significant progress in birational geometry in positive characteristic in dimension $\leq 3$, such as \cite{HaconXu}. In recent years,  input from $p$-adic Hodge theory has made it possible to prove similar vanishing theorems in mixed characteristic algebraic geometry; this has led to solutions of longstanding questions in commutative algebra and also to progress in the minimal model program in mixed characteristic.

\subsection{Vanishing theorems in commutative algebra}
\label{ss:VanThm}

$F$-singularity theory is the study of singularities in positive characteristic via the behaviour of the Frobenius endomorphism. It was born with a classical theorem of Kunz \cite{KunzRegular} proving that a noetherian $\mathbf{F}_p$-algebra is regular exactly when its Frobenius endomorphism is flat.  This subject was systematically developed by Hochster--Huneke and several others over many decades; see \cite{HHTightClosure} as well as the survey \cite{TakagiWatanabeFSingSurvey}. An important landmark  in the subject was a Cohen--Macaulayness result of Hochster--Huneke \cite{HHBigCM}; see \cite{HunekeAICSurvey} for a fairly recent survey. The following recent result extends this to mixed characteristic:

\begin{theorem}[Cohen--Macaulayness of $R^+$, \cite{BhattCM, BMPSTWW}]
\label{RPlusCM}
 Let $R$ be an excellent noetherian domain. Let $R^+$ be the integral closure of $R$ in an algebraic closure of its fraction field. Then the $p$-adic completion $\widehat{R^+}$ is Cohen--Macaulay over $R$. 
 \end{theorem}

\begin{remark}[A concrete formulation]
\label{rmk:ElementaryCM}
Despite involving the large ring $R^+$, Theorem~\ref{RPlusCM} is a finitistic statement whose essential content is the following: if $R$ is local and $\underline{x} := \{p,x_1,...,x_d\}$ is a system of parameters, then any relation on $\underline{x}$ becomes a linear combination of the trivial Koszul relations in a finite extension $S$ of $R$. This formulation explains why Theorem~\ref{RPlusCM} can be viewed as a ``vanishing theorem up to finite covers'': it says that the local cohomology classes on $R$ coming from the potentially non-trivial relations can be annihilated by passing to finite extensions $R \to S$. Moreover, it also highlights the essential difficulty: one must {\em construct} finite extensions of $R$ from the unwanted relations. 
\end{remark}

\begin{remark}[Weakly functorial Cohen--Macaulay algebras]
\label{rmk:BigCM}
Andr\'e's recent resolution  \cite{AndreDSC,AndreAbhyankar} of Hochster's direct summand conjecture led to a lot of activity in mixed characteristic commutative algebra, including \cite{BhattDSC,MaSchwedeBoundSymbolic,MaSchwedeBCM,HeitmannMaTor,ShimomotoBCM,MSTWWAdjoint}; see \cite{MaSchwedeSurvey} for a recent survey. In particular, Andr\'e \cite{AndreWeakFun} and Gabber \cite{GabberBCM} proved the existence of ``weakly functorial Cohen--Macaulay algebras'' in the key remaining mixed characteristic case (via rather indirect constructions). This existence result implies many of the ``homological conjectures'' in commutative algebra (a notable exception being Serre's intersection multiplicity conjecture); see \cite{HochsterHomological,AndreICM,HochsterHomologicalMSRI}.
 Prior to Andr\'e and Gabber's work, this existence was  known (\cite{HochsterHeitmannBigCM,HeitmannDSC3}) only in dimension $\leq 3$. Theorem~\ref{RPlusCM} now yields an alternative and extremely simple construction of such weakly functorial Cohen--Macaulay algebras in mixed characteristic:  we may simply use $\widehat{R^+}$. 
\end{remark}

\begin{remark}[What was known?]
Theorem~\ref{RPlusCM} is straightforward in dimension $\leq 2$, and is the main result of \cite{HHBigCM} in the positive characteristic case.  In mixed characteristic, Theorem~\ref{RPlusCM} is new even in dimension $3$: it was previously known \cite{HeitmannDSC3} in dimension $\leq 3$ only in the almost category  (in the sense of Faltings' almost mathematics \cite{FaltingspHT,GabberRamero}); see \cite[Remark 1.9]{BhattCM} for an explanation of prior expectations. 
\end{remark}

\begin{remark}[Splinters]
A noetherian commutative ring $R$ is called a {\em splinter} if it satisfies the conclusion of the direct summand conjecture, i.e., it splits off as a module from every finite extension. This class of singularities, formally introduced in \cite{MaSplinter}, has recently received renewed attention (e.g., \cite{MaVanishingTor,AntieauDattaValuation,DattaTuckerOpen,AndreFiorot}). An external reason to care about this notion is a major conjecture in $F$-singularity theory (\cite[page 85]{HHBigCM}, \cite[page 640]{HHTightClosureJAG}): splinters in characteristic $p$ are expected to be the same as strongly $F$-regular rings (see \cite[end of \S 3]{MaPolstraBook} for a discussion). This conjecture is known for $\mathbf{Q}$-Gorenstein rings \cite{SinghQGorSplinter}. One consequence of this conjecture is that characteristic $p$ splinters are  derived splinters, i.e., they satisfy a derived version of the splinter condition for any proper surjective map and are thus analogous to rational singularities. This consequence was proven unconditionally in \cite{BhattDDSCPosChar}. Methods from \cite{BhattLurieRH1} used in proving Theorem~\ref{RPlusCM} give the same result in mixed characteristic. In conjunction with Theorem~\ref{RPlusCM} itself, one learns that that any mixed characteristic splinter is Cohen--Macaulay and has rational singularities in the sense of \cite{KovacsRational}; it would be interesting to prove the latter (even just after inverting $p$) without using $p$-adic Hodge theory.
\end{remark}

\begin{remark}[Ingredients in the proof of Theorem~\ref{RPlusCM}]
Using essentially elementary methods, \cite{BMPSTWW} reduces Theorem~\ref{RPlusCM} to the statement that $R^+/p$ is Cohen--Macaulay over $R/p$, which is proven in \cite{BhattCM}. Despite the simple reformulation highlighted in Remark~\ref{rmk:ElementaryCM}, the proof relies on two major theoretical inputs. The first is prismatic cohomology (Theorem~\ref{PrismaticCoh}), which gives a substitute for the Frobenius operator in mixed characteristic; this allows one to begin mimicking the cohomological proof of \cite{HHBigCM} given in \cite{HunekeLyubeznik} in mixed characteristic at the cost of replacing rings with derived rings. The second is the $p$-adic Riemann--Hilbert functor from Theorem~\ref{RHtorsion} below, applied to certain perverse $\mathbf{F}_p$-sheaves on the generic fibre $\mathrm{Spec}(R[1/p])$ arising from finite covers, to facilitate the induction on dimension strategy of \cite{HunekeLyubeznik}. This proof is not effective, and it might be interesting to explicitly constructed the relevant covers in low dimensional examples, such as cones over smooth projective curves and surfaces over a $p$-adic discrete valuation ring.
\end{remark}

\subsection{Birational geometry}
\label{ss:BG}

There is a well-known analogy between projective geometry and local algebra, e.g., the global cohomological properties of a projective variety $X \subset \mathbf{P}^n$ are faithfully reflected  in the local cohomological properties of its affine cone $Y \subset \mathbf{A}^{n+1}$ over $X$ near the vertex $0 \in Y$. This analogy suggests that Theorem~\ref{RPlusCM} ought to have a global variant; this is indeed the case, and the result can be summarized as follows:

\begin{theorem}[Kodaira vanishing up to finite covers, \cite{BhattCM}]
\label{thm:KV}
Let $V$ be a $p$-adic discrete  valuation ring (e.g., $V=\mathbf{Z}_p$). Let $X/V$ be a flat proper scheme equipped with a semiample and big line bundle $L$. Then any $p$-power torsion class in $H^*(X,L^{-1})$, $H^*(X,\mathcal{O}_X)$, or $H^*(X,L)$ can be annihilated by pullback to a finite cover of $X$.
\end{theorem}

Analogous results hold true in the relative setting \cite{BhattCM}, and were previously known in characteristic $p$ (\cite{HHBigCM} for $L$ ample, and \cite{BhattDDSCPosChar} in general). 

\begin{remark}[Relation to Kodaira vanishing]
The classical Kodaira vanishing theorem says that 
\[ H^{< \dim(Y)}(Y,M^{-1}) = 0\]
 for a smooth projective variety $Y/\mathbf{C}$ with ample line bundle $M$. This assertion is false in characteristic $p$ (\cite{RaynaudKodaira}) and mixed characteristic (by Totaro, see \cite[Footnote 1]{BMPSTWW}). The $L^{-1}$ case of Theorem~\ref{thm:KV} can be viewed as an ``up to finite covers'' variant of the Kodaira vanishing theorem that {\em is} true in mixed characteristic: spurious cohomology classes --- ones that should not be there if Kodaira vanishing were true for $(X,L)$ --- can be annihilated by passing to finite covers. This ``up to finite covers'' perspective was pioneered in characteristic $p$ by \cite{SmithFiniteCover} in the wake of \cite{HHBigCM}. 

For completeness, we remark that an ``up to finite covers'' version of the more general Kodaira--Akizuki--Nakano vanishing theorem also holds true in the setting of Theorem~\ref{thm:KV}: in fact, the cases not covered by Theorem~\ref{thm:KV} are much easier as sheaves of differential forms themselves become $p$-divisible on passage to finite covers. 
\end{remark}

\begin{remark}[Relation to the $p$-adic Poincar\'e lemma]
The assertion in Theorem~\ref{thm:KV} for $H^*(X,\mathcal{O}_X)$,  with finite covers weakened to alterations, was previously known by \cite{Beilinsonpadic,BhattMixedCharPDiv}; in fact, it formed the key geometric ingredient in the proof of the $p$-adic Poincar\'e lemma in \cite{Beilinsonpadic}. Curiously, while the $p$-adic Poincar\'e lemma was used in \cite{Beilinsonpadic} to give a new proof of the fundamental de Rham comparison conjecture in $p$-adic Hodge theory, the proof of Theorem~\ref{thm:KV} uses the full strength of modern advances in $p$-adic Hodge theory (such as the primitive comparison theorem of \cite{ScholzepHT} for arbitrarily singular varieties). 
\end{remark}

We end this section with an application of Theorem~\ref{thm:KV} to birational geometry in mixed characteristic. Briefly, it is possible to use this variant of Kodaira vanishing in a critical lifting argument in an inductive proof of the existence of flips in dimension $3$, following  \cite{HaconXu} (which goes back to ideas of Shokurov). Combining this with Witaszek's recent mixed characteristic  analog \cite{WitaszekKeel} of Keel's semiampleness theorems \cite{KeelSemiample},  it  became possible to emulate the ideas of \cite{HaconXu,Birkar16,BW17,DW19}  (amongst others) to show the following:

\begin{theorem}[Minimal model program in mixed characteristic, \cite{BMPSTWW, TakamatsuYoshikawaMMP}]
\label{MMP}
One can run the minimal model program for arithmetic threefolds whose residue characteristics are $> 5$.
\end{theorem}

Theorem~\ref{MMP} uses ideas from \cite{BSTAlterations, MSTWWAdjoint} and extends \cite{TanakaSurface,KawamataMMP}. Global geometric applications of (the ideas going into) Theorem~\ref{MMP} can be found in \cite{BMPSTWW, TakamatsuYoshikawaMMP,HaconLamarcheSchwede,StigantMori,XieXueMMP4}.

\begin{remark}[The $+$-stable sections]
\label{rmk:B0Localize}
We informally discuss a new notion introduced in the proof of Theorem~\ref{MMP} in a simple case, and state a question; see \cite[\S 4]{BMPSTWW} or \cite[\S 3.2]{TakamatsuYoshikawaMMP} for the general notion. For reasonable mixed characteristic rings $R$, one can define a submodule $B^0(R,\omega_R) \subset \omega_R$ of the dualizing module $\omega_R$: it is the submodule of elements that lift to all alterations of $\mathrm{Spec}(R)$ under the trace maps. If $R$ is regular, then $B^0(R,\omega_R) = \omega_R$, so in general $B^0(R,\omega_R)$ is an invariant measuring the singularities of $R$. Analogous invariants exist in characteristic $0$ (given by the Grauert--Riemenschneider sheaf \cite[Example 4.3.12]{LazarsfeldPos1}) and characteristic $p$ (given by the parameter test submodule \cite[\S 2.5 \& Corollary 3.4]{BSTAlterations}). Basic properties of $B^0(R,\omega_R)$, such as its behaviour under alterations or restriction to divisors, play a key role in the proof of Theorem~\ref{MMP}. However, a fundamental question about these invariants remains open: does their formation commute with localization? Due to the infinite intersection implicit in the definition of $B^0(R,\omega_R)$, this question is  delicate. Nevertheless, a positive answer (which we expect) would have several geometric applications. As evidence for a positive answer, using Theorem~\ref{RH2}, one can show the claim for inverting $p$: the localization $B^0(R,\omega_R)[1/p]$ agrees with the Grauert--Riemenschneider sheaf of $\mathrm{Spec}(R[1/p])$ (work in progress as a sequel to \cite{BMPSTWW}). We refer to \cite[\S 8]{HaconLamarcheSchwede} for more discussion of this question.
\end{remark}

\section{$p$-adic Riemann--Hilbert}
\label{sec:RH}

The Riemann--Hilbert problem has a rich history, going back at least to Hilbert's $21$st problem. In modern terms, it asked if any $\mathbf{C}$-local system on a smooth complex algebraic curve $X$ can  be realized as the solution system of a flat vector bundle on $X$ with regular singularities at $\infty$; this variant was (precisely formulated and) solved by Deligne in \cite{DeligneRH}, which also contained a higher dimensional analog.  Soon after, this picture was significantly generalized  by Kashiwara and Mebkhout: there is a natural equivalence of categories between  topological objects ($\mathbf{C}$-linear perverse sheaves) and differential objects (regular holonomic $\mathcal{D}$-modules) on any smooth complex variety, see \cite{BorelDMod}. 

In this section, we discuss joint work with Lurie towards a $p$-adic analog of the preceding story; our aim was to extend existing results attaching flat connections to $p$-adic local systems on $p$-adic varieties  (such as \cite{FaltingsSimpson1,FaltingsSimpson2,AbbesGrosTsujiSimpson,ScholzepHT,LiuZhu,DLLZRH})  to  $p$-adic constructible complexes and in particular, to $p$-adic perverse sheaves. Unlike the complex picture, there are several meanings one can attach to ``$p$-adic sheaves'': one can work with $\mathbf{Z}/p^n$, $\mathbf{Z}_p$ or $\mathbf{Q}_p$-coefficients. Our theorem for $\mathbf{F}_p$-coefficients is the following (the $\mathbf{Z}/p^n$ case is analogous):

\begin{theorem}[Riemann--Hilbert for torsion coefficients, \cite{BhattLurieRH1}]
\label{RHtorsion}
Let $C/\mathbf{Q}_p$ be a complete and algebraically closed extension. Let $X/\mathcal{O}_C$ be a finite type scheme. Then there is a natural exact functor
\[ \mathrm{RH}:D^b_{\mathrm{cons}}(X_C, \mathbf{F}_p) \to D^b_{\mathrm{qc}}(X \otimes_{\mathcal{O}_C} \mathcal{O}_C/p). \]
This functor commutes with proper pushforward, intertwines Verdier and Grothendieck duality in the almost category \cite{FaltingspHT,GabberRamero}, and interacts well with the perverse $t$-structure. 
\end{theorem}

The functor $\mathrm{RH}$ above also almost commutes with tensor products and pullbacks provided the target is refined to $\mathrm{RH}(\mathbf{F}_p)$-modules. In fact, it is possible to refine the target further to Frobenius modules over the tilt $\mathrm{RH}(\mathbf{F}_p)^\flat$; the resulting functor is fully faithful, and agrees with the construction in \cite{BhattLurieModpRH} (which was a dual form of \cite{EmertonKisin} that works for all characteristic $p$ schemes) when $X$ has characteristic $p$.

\begin{remark}[Relation to existing work in $p$-adic geometry]
Theorem~\ref{RHtorsion} appears to be the first general construction attaching coherent objects to constructible $\mathbf{F}_p$-sheaves on algebraic varieties in characteristic $0$.  On the other hand, several ingredients that go into the proof have appeared before in $p$-adic arithmetic geometry. Indeed, the functor $\mathrm{RH}$ can be regarded as a generalization of a perfectoidization functor from Remark~\ref{rmk:Perfections} to non-constant coefficients: one can almost identify $\mathrm{RH}(\mathbf{F}_p)$ with $\mathcal{O}_{X,\mathrm{perfd}}/p$. Moreover, the compatibility with duality with constant coefficients is closely related to the Gabber--Zavyalov approach \cite{ZavyalovPD} to Poincar\'e duality for the $\mathbf{F}_p$-cohomology of rigid spaces. Nevertheless, the flexibility of applying $\mathrm{RH}(-)$ to non-constant perverse coefficients is immensely useful in applications including Theorems~\ref{RPlusCM} and \ref{thm:KV} or the localization result mentioned in Remark~\ref{rmk:B0Localize}. Relatedly, let us mention that Theorem~\ref{RHtorsion} itself suffices to prove Theorem~\ref{RPlusCM} in the almost category, extending Heitmann's almost vanishing theorem \cite{HeitmannDSC3} to arbitrary dimensions.
\end{remark}

Prima facie, Theorem~\ref{RHtorsion} looks quite different from the complex Riemann--Hilbert correspondence: the output is a quasicoherent (and in fact almost coherent) complex rather than a $\mathcal{D}$-module. In fact, the functor in Theorem~\ref{RHtorsion} is better understood as a $p$-adic analog of a construction from Saito's fundamental work \cite{SaitoMHM} on mixed Hodge modules. Recall that this theory gives a {\em filtered} refinement of the classical Riemann--Hilbert functor for many constructible sheaves, including those that are ``of geometric origin''. More precisely, given a smooth proper complex variety $X$, any mixed Hodge module on $X$ has an underlying $\mathcal{D}_X$-module equipped with a Hodge filtration as well as an underlying perverse sheaf; the picture relating them can be summarized in the following commutative diagram:
\[\xymatrix{ 
& \mathrm{MHM}(X) \ar[ld]_-{\text{forget}} \ar[rd]^-{\text{forget}} & & & \\ 
 D^b_{\text{cons}}(X,\mathbf{C}) \ar[rd]^-{\mathrm{RH}^{\text{cl}}} & &  DF_{\mathrm{coh}}(\mathcal{D}_X) \ar[r]^-{\mathrm{gr}^*(-)} \ar[ld]_-{\text{forget}} & D^b_{\mathrm{coh},\mathrm{gr}}(T^* X) \ar[r]^-{\Omega^*(-)} & D^b_{\mathrm{coh},\mathrm{gr}}(X) \\
& D^b(\mathcal{D}_X),  }\]
where $\mathrm{MHM}(X)$ is Saito's category of mixed Hodge modules, $DF_{\mathrm{coh}}(\mathcal{D}_X)$ is a suitable derived category of $\mathcal{D}_X$-modules equipped with a ``good'' filtration, the functor $\mathrm{RH}^{\text{cl}}$ is the classical Riemann--Hilbert functor, the functor $\mathrm{gr}^*(-)$ is the associated graded construction carrying a filtered $\mathcal{D}$-module to a graded $\mathcal{O}_X$-module with an action of $\mathrm{gr}^* \mathcal{D}_X = \mathrm{Sym}^* (T_X)$ (i.e., a Higgs module), and the functor $\Omega^*(-)$ is the graded Higgs complex construction. Heuristically, the functor in Theorem~\ref{RHtorsion} is an analog of the composite correspondence 
\begin{equation}
\tag{$\widetilde{\mathrm{RH}}$}
\label{Saito}
 D^b_{\text{cons}}(X, \mathbf{C}) \xleftarrow{\text{forget}} \mathrm{MHM}(X) \xrightarrow{\Omega^*(-) \circ \mathrm{gr}^*(-) \circ \text{forget}} D^b_{\mathrm{coh},\mathrm{gr}}(X)
 \end{equation}
for $\mathbf{F}_p$-coefficients. Slightly surprisingly, unlike in the complex story, we get an honest functor instead of a correspondence in the $p$-adic setting. (On the other hand, objects of $\mathrm{MHM}(X)$ also have a weight filtration, which we ignore in our discussion.)

\begin{remark}[Why is there no grading?]
In comparison with the correspondence \eqref{Saito}, there is no grading in the target of Theorem~\ref{RHtorsion}. But this is to be expected: the grading on the target of \eqref{Saito} reflects the fact that objects in $\mathrm{MHM}(X)$ are fairly motivic in nature, e.g., they give variations of Hodge structures on an open subset of $X$. In contrast, in Theorem~\ref{RHtorsion} we are working with {\em all} constructible sheaves over the algebraically closed field $C$, so there is no motivicity or even a Galois action.
\end{remark}

The previous discussion suggests it might be useful to lift Theorem~\ref{RHtorsion} to $\mathbf{Q}_p$-coefficients and restrict to sheaves defined over a discretely valued field (so there is a Galois action) in order to obtain a $p$-adic variant of \eqref{Saito}. This can indeed be done, and the resulting structure seems slightly cleaner than \eqref{Saito}:

\begin{theorem}[Riemann--Hilbert for $\mathbf{Q}_p$-coefficients, \cite{BhattLurieRH2}]
\label{RH2}
Let $K/\mathbf{Q}_p$ be a finite extension. Let $X/K$ be a smooth proper variety. Then there is a natural exact functor
\[ \mathrm{RH}_{\mathcal{D}}:D^b_{\mathrm{wHT}}(X, \mathbf{Q}_p) \to DF_{\mathrm{coh}}(\mathcal{D}_X), \]
where the source is a full subcategory of $D^b_{\mathrm{cons}}(X, \mathbf{Q}_p)$ spanned by what we call ``weakly Hodge--Tate sheaves'' (including all sheaves of geometric origin). This functor commutes with proper pushforward, intertwines Verdier and Grothendieck duality, and interacts well with the perverse $t$-structure.
\end{theorem}

Theorem~\ref{RH2} represents ongoing work in progress with Lurie, and the statement above is not quite optimal (e.g., there is a variant for singular $X$).

\begin{remark}[The case of local systems]
The functor $\mathrm{RH}_{\mathcal{D}}$ from Theorem~\ref{RH2} is not really new for local systems: up to a certain nilpotent operator encoding that a weakly Hodge--Tate local system is not quite de Rham, it coincides with the one appearing in \cite[Theorem 1.5]{LiuZhu} (and is thus related to constructions from \cite{ScholzepHT}; see also \cite{FaltingsSimpson1,AbbesGrosTsujiSimpson,DLLZRH}). However, for geometric applications such as Example~\ref{ex:GRVan} below, it is critical to apply $\mathrm{RH}_{\mathcal{D}}$ to constructible complexes that are not local systems.
\end{remark}

\begin{remark}[Why is the Hodge filtration automatic?]
Theorem~\ref{RH2} implies that constructible $\mathbf{Q}_p$-sheaves $F$ of geometric origin on a variety $X/K$ as above have a functorially attached filtered $\mathcal{D}_X$-module $\mathcal{M} := \mathrm{RH}_{\mathcal{D}}(F)$, i.e., the Hodge filtration on the $\mathcal{D}_X$-module $\mathcal{M}$ is actually determined by $F$, unlike in the correspondence \eqref{Saito}. This discrepancy is ultimately because the constructible sheaves in Theorem~\ref{RH2} carry Galois symmetries as they are defined over $K$. Moreover, this is perfectly consistent with known phenomena in $p$-adic Hodge theory that stem ultimately from the richness of the absolute Galois group $G_K$ of $K$. For instance, when $X=\mathrm{Spec}(K)$ and $F = Rf_* \mathbf{Q}_p$ for a smooth proper map $f:Y \to X$, we are simply observing that the $G_K$-representation $H^*(Y_{\overline{K}}, \mathbf{Q}_p)$ knows the de Rham cohomology of $H^*_{\mathrm{dR}}(Y/K)$ as a filtered vector space (and in particular knows the Hodge numbers of $X$) via the de Rham comparison; see \cite{ItoBirInv} for a purely geometric application of this fact.
\end{remark}

As Theorem~\ref{RH2} gives an honest functor, one can now directly apply $\mathrm{RH}_{\mathcal{D}}$ to deep theorems on the constructible side, such as the BBDG decomposition theorem \cite{BBDG}, to obtain highly nontrivial results on the coherent side. This mechanism appears robust enough to yield some results in birational geometry that are traditionally best understood via mixed Hodge module theory, e.g.,  Koll\'ar's vanishing theorems \cite{KollarHigher1,KollarHigher2} (see \cite[\S 25]{SchnellMHMOverview} for the Hodge module proof); we sketch the argument for vanishing next to illustrate this idea.

\begin{example}[Recovering Koll\'ar vanishing, $p$-adically]
\label{ex:GRVan}
Fix a finite extension $K/\mathbf{Q}_p$. Say $f:Y \to X$ is a projective surjective morphism of proper $K$-varieties of dimensions $d_Y$ and $d_X$ respectively with $Y$ smooth. Consider the functor
\[ \mathrm{RH}:D^b_{\mathrm{wHT}}(Y, \mathbf{Q}_p) \to D^b_{\mathrm{coh},\mathrm{gr}}(Y)\]
obtained by composing the functor $\mathrm{RH}_{\mathcal{D}}$ from Theorem~\ref{RH2} with $\Omega^*(-) \circ \mathrm{gr}^*(-)$, as in \eqref{Saito}. Essentially by the local Hodge--Tate decomposition of \cite{ScholzepHT}, we have 
\[ \mathrm{RH}(\mathbf{Q}_p[d_Y]) = \bigoplus_i \Omega^i_{Y/K}[d_Y-i]\] 
with its natural grading, so $i$-forms have weight $i$. (If $Y$ were singular, one would have a similar formula with the Deligne--Du Bois variants $\underline{\Omega}^i_{Y/K}$ of differential forms, as in \cite{DuBoisFilt} and  \cite[\S 7.3]{PetersSteenbrink}, on the right by \cite{GuoHodgeTate}.) Pushing forward along $f$, using the proper pushforward compatibility of $\mathrm{RH}$, and extracting the weight $d_Y$ summand gives
\[\mathrm{RH}(Rf_* \mathbf{Q}_p[d_Y])_{\mathrm{wt}=d_Y} = Rf_* \omega_Y.\]
On the other hand, the decomposition theorem \cite{BBDG, DeligneDecompDerived} shows that 
\[ Rf_* \mathbf{Q}_p[d_Y] \simeq \left(\bigoplus_{i=-(d_Y-d_X)}^{d_Y-d_X} {}^p \mathcal{H}^i[-i]\right) \oplus N\]
where each ${}^p \mathcal{H}^i$ is perverse and $N$ is a summand of $Rg_* \mathbf{Q}_p[d_Y]$ with $g:Y_Z \to Z \subset X$ being the restriction of $f$ over the closed subvariety $Z \subsetneq X$ where $f$ is not smooth. The singular variant of the reasoning just used for $f$ applied to $g$ then shows that
\[ \mathrm{RH}(Rg_* \mathbf{Q}_p[d_Y])_{\mathrm{wt}=d_Y} = Rg_* \underline{\Omega}^{d_Y}_{Y_Z}  \simeq 0,\]
where the last vanishing follows as $\underline{\Omega}^{d_Y}_{Y_Z} = 0$ since $d_Y > \dim(Y_Z)$ (see \cite{GuoHodgeTate} for a purely $p$-adic proof of this property of Deligne-Du Bois complexes). But then the same vanishing is also true for the summand $N$ of $Rg_* \mathbf{Q}_p[d_Y]$, so we learn that 
\[ Rf_* \omega_Y = \mathrm{RH}(Rf_* \mathbf{Q}_p[d_Y])_{\mathrm{wt}=d_Y} = \bigoplus_{i=-(d_Y-d_X)}^{d_Y-d_X} \mathrm{RH}(\mathcal{H}^i[-i])_{\mathrm{wt}=d_Y}.\]
The perverse exactness properties of $\mathrm{RH}$ now imply that the $i$-th summand on the right lies in $D^{\leq i}$ whence $Rf_* \omega_Y \in D^{\leq d_Y-d_X}$ as $i \leq d_Y-d_X$, i.e., 
\[ R^j f_* \omega_Y = 0 \quad \text{for} \quad j > d_Y-d_X,\]
proving the Koll\'ar vanishing theorem \cite[Theorem 2.1]{KollarHigher1}. From this perspective, one answer to Koll\'ar's question ``Why is $\omega_Y$ better behaved than $\mathcal{O}_Y$?" \cite{KollarHigher2} could be the following: as $\omega_Y$ is the highest Hodge--Tate weight summand  of $\mathrm{RH}(\mathbf{Q}_p[d_Y])$, it does not see interference from smaller dimensional varieties when moved around via operations such as $Rf_*$.
\end{example}

\subsection*{Acknowledgements}
The landscape of algebraic geometry in mixed characteristic has changed substantially in the last decade. I feel extremely fortunate to have had the opportunity to witness many aspects of this metamorphosis up close. I am thus eternally indebted to my mentors, collaborators, and friends for their support, wisdom and generosity (with time and ideas) over the years.   Thanks also to Johan de Jong, H\'el\`ene Esnault, Valia Gazaki, Toby Gee, Lars Hesselholt, Jacob Lurie, Linquan Ma, Akhil Mathew, Davesh Maulik, Mircea Mustata, Arthur Ogus, Alex Perry, Peter Scholze, Karl Schwede, Kevin Tucker, and especially Wei Ho, Luc Illusie, Matthew Morrow and Anurag Singh for numerous helpful comments during the preparation of this survey.

This work was partially supported by the NSF (\#1801689, \#1952399, \#1840234), the Packard Foundation, and the Simons Foundation (\#622511).


\bibliographystyle{abbrv}
\bibliography{mybib}

\end{document}